\documentclass{article}

\usepackage{latexsym,amsfonts,amsmath,amssymb,graphicx,hyperref,verbatim}

\newcommand{\EE}{\mbox{\bf E}\,}
\newcommand{\PP}{\mbox{\bf P}\,}
\newcommand{\QQ}{\mbox{\bf Q}\,}

\newcommand{\R}{\mathbb{R}}
\newcommand{\C}{\mathbb{C}}
\newcommand{\Q}{\mathbb{Q}}
\newcommand{\HH}{\mathbb{H}}
\newcommand{\N}{\mathbb{N}}

\newcommand{\Z}{\mathbb{Z}}

\newcommand{\pa}{\partial}

\newcommand{\F}{{\cal F}}
\newcommand{\no}{\noindent}

\newcommand{\BGE}{\begin{equation}}
\newcommand{\BGEN}{\begin{equation*}}
\newcommand{\EDE}{\end{equation}}
\newcommand{\EDEN}{\end{equation*}}

\def\eps{\varepsilon}
\def\til{\widetilde}
\def\ha{\widehat}
\def\sem{\setminus}
\def\Om{\Omega}
\def\lin{\overline}
\def\vphi{\varphi}

\def\del{\delta}

\def\h0{{\bf h}}

\DeclareMathOperator{\hcap}{hcap} \DeclareMathOperator{\id}{id}
\DeclareMathOperator{\Imm}{Im } \DeclareMathOperator{\Ree}{Re }

 \DeclareMathOperator{\HP}{HP}

\newtheorem{Lemma}{Lemma}[section]
\newtheorem{Theorem}{Theorem}[section]
\newtheorem{Definition}{Definition}[section]
\newtheorem{Corollary}{Corollary}[section]
\newtheorem{Proposition}{Proposition}[section]
\newtheorem{Conjecture}{Conjecture}
\numberwithin{equation}{section}

\begin{document}

\title{\bf Reversibility of Some Chordal SLE$(\kappa;\rho)$ Traces}
\date{\today}
\author{Dapeng Zhan\footnote{Yale University}}
\maketitle
\begin{abstract} We prove that, for $\kappa\in(0,4)$ and $\rho\ge (\kappa-4)/2$,
 the chordal SLE$(\kappa;\rho)$ trace started from $(0;0^+)$ or $(0;0^-)$ satisfies
 the reversibility property. And we obtain the equation for the reversal of
 the  chordal SLE$(\kappa;\rho)$ trace started from $(0;b_0)$, where $b_0>0$.
\end{abstract}

\section{Introduction}

In the proof of the reversibility of the SLE$(\kappa)$ trace (\cite{reversibility}), where $\kappa\in(0,4]$,
a new technique was developed to construct a coupling of two SLE$(\kappa)$ traces,
such that in that coupling, the images of the two traces coincide, and the directions of the two traces are opposite.
That technique was then used to prove the Duplantier's duality conjecture (\cite{duality}\cite{duality2}). Comparing
Theorem 5.4 in \cite{duality} with Julien Dub\'edat's Conjecture 2 in \cite{Julien-Duality}, the author proposed the following
conjecture in \cite{duality}.

\begin{Conjecture} Let $\beta_0(t)$, $0\le t<\infty$, be a  chordal SLE$(\kappa; \rho_+, \rho_-)$
 trace started from $(0; 0^+, 0^-)$, where $\kappa\in(0,4)$ and $\rho_+,\rho_-\ge(\kappa-4)/2$.
 Let $W_0(z)=1/\lin{z}$. Then after a time-change, $(W_0(\beta_0(1/t)))$, $0<t<\infty$, has the same distribution as $(\beta_0(t))$, $0<t<\infty$.\label{conjec}
\end{Conjecture}

It's already known that this conjecture holds in some special cases. If $\rho_+=\rho_-=0$, then $\beta_0$ is a
standard SLE$(\kappa)$ trace, and the result follows from \cite{reversibility}.
If $\kappa=0$, then $\beta_0$ is a half line from $0$ to $\infty$, which is a
trivial case. If $\kappa=4$, then it follows from the convergence of the discrete Gaussian free field contour line
(\cite{SS}); and it is also a special case of Theorem 5.5 in \cite{duality}. The motivation of the current paper
is to prove the above conjecture. We will only prove part of it, that is, the case when $\rho_+$ or $\rho_-$ equals to $0$.
If, for example, $\rho_-=0$, then $\beta_0$ reduces to a chordal SLE$(\kappa; \rho_+)$ trace started from $(0; 0^+)$.
The main theorem of this paper is the following.

\begin{Theorem} Let $\kappa\in(0,4)$ and $\rho\ge (\kappa-4)/2$. Suppose $\beta_0(t)$, $0\le t<\infty$, is a
chordal SLE$(\kappa;\rho)$ trace started from $(0;0^\sigma)$, where $\sigma\in\{+,-\}$. Let $W_0(z)=1/{\lin z}$. Then after a time-change,
$W_0(\beta_0(1/t))$, $0<t<\infty$, has the same distribution as $\beta_0(t)$, $0<t<\infty$.
\label{reversal*}\end{Theorem}

We will see that Theorem \ref{reversal*} here and Theorem 5.4 in \cite{duality} imply Dub\'edat's conjecture. Besides the special cases
that $\rho=0$, $\kappa=0$ or $4$, the above theorem is also known to be true in the case that $\kappa=8/3$. This follows from
\cite{LSW-8/3} because the image of $\beta_0$ satisfies the left-sided or right-sided restriction property with exponent
depending on $\rho$, and the one-sided restriction measure is invariant under the map $W_0(z)=1/{\lin z}$.

The proof of Theorem \ref{reversal*} will be completed in the last section.
We will use the technique used in \cite{reversibility} and \cite{duality}.
The new difficulty here is that when applying  the above technique, we need some information about the ``middle'' part
of the curve $\beta_0$. This means that given a stopping time $T_1>0$ and a ``backward'' stopping time $T_2<\infty$ with $T_1<T_2$,
we need to know the conditional distribution of $\beta_0(t)$, $T_1< t < T_2$, given the curves $\beta_0((0; T_1])$ and
 $\beta_0([T_2;\infty))$. This is known in some special cases. If $\beta_0$ is a standard chordal SLE$(\kappa)$ trace, which
 corresponds to the case that $\rho=0$, then $\beta_0(t)$, $T_1< t < T_2$, is a time-change of a chordal SLE$(\kappa)$ trace
 in $\HH\sem(\beta_0((0; T_1])\cup \beta_0([T_2;\infty)))$ from $\beta_0(T_1)$ to $\beta_0(T_2)$. If $\kappa=4$, from the
 proof of Theorem 5.5 in \cite{duality}, we see that $\beta_0(t)$, $T_1< t < T_2$, is a time-change of a generic
 SLE$(\kappa;\rho)$ trace in $\HH\sem(\beta_0((0; T_1])\cup \beta_0([T_2;\infty)))$. In the general case, as we will see,
 the conditional distribution of $\beta_0(t)$, $T_1< t < T_2$, is complicated. To describe this middle part of $\beta_0$,
 we will use hypergeometric functions to define a new kind of SLE-type processes, which are called intermediate SLE$(\kappa;\rho)$
 processes. These new SLE-type processes will also be used to describe the reversal of an SLE$(\kappa;\rho)$ trace
 whose force point is not degenerate. This is Theorem \ref{reversal-2} below, whose proof will also be completed in the last section.

\begin{Theorem} Suppose $\beta_0(t)$, $0\le t<\infty$, is a  chordal SLE$(\kappa;\rho)$ trace started from $(0;b_0)$ with $b_0>0$.
Let $W_0(z)=1/{\lin z}$. Then after a time-change, $W_0(\beta_0(1/t))$, $0<t<\infty$, has the same distribution as
a degenerate  intermediate SLE$(\kappa;\rho)$ trace with force points $0^+$ and $1/b_0$.
\label{reversal-2}\end{Theorem}

The current paper will frequently use results from \cite{reversibility} and \cite{duality}. The reader is suggested to have copies
of those two papers by hand for convenience.

After finishing the first version of this paper, the author noticed that Corollary 9 in \cite{Julien-Duality-2} is equivalent
to Theorem \ref{reversal*} here. It seems to the author that some important details are omitted in \cite{Julien-Duality-2}. 
The proofs in this paper will be completed, and contain all details. And the approach of this paper is somewhat different from 
that in \cite{Julien-Duality-2}.

\section{Preliminary} If $H$ is a bounded and relatively  closed
subset of $\HH=\{z\in\C:\Imm z>0\}$, and $\HH\sem H$ is simply
connected, then we call $H$ a hull in $\HH$ w.r.t.\ $\infty$. For
such $H$, there is $\vphi_H$ that maps $\HH\sem H$ conformally onto
$\HH$, and satisfies $\vphi_H(z)=z+\frac{c}{z}+O(\frac 1{z^2})$ as
$z\to\infty$, where $c=\hcap(H)\ge 0$ is called the half-plane  capacity of $H$.
A hull $H$ with $\hcap(H)=c$ has diameter at least $\sqrt c$.
If $H_1\subset H_2$ are hulls
in $\HH$ w.r.t.\ $\infty$, then $H_2/H_1:=\vphi_{H_1}(H_2\sem H_1)$
is also a hull in $\HH$ w.r.t.\ $\infty$, and we have
$\vphi_{H_2}=\vphi_{H_2/H_2}\circ \vphi_{H_1}$.

For a real interval $I$, we use $C(I)$ to denote the space of real
continuous functions on $I$. For $T>0$ and $\xi\in C([0,T))$, the
chordal Loewner equation driven by $\xi$ is
\BGE \pa_t\vphi(t,z)=\frac{2}{\vphi(t,z)-\xi(t)},\quad \vphi(0,z)=z.\label{chordal-equation}\EDE
For $0\le t<T$, let $K(t)$ be the set of $z\in\HH$ such that the
solution $\vphi(s,z)$ blows up before or at time $t$. Then each $K(t)$
is a hull in $\HH$ w.r.t.\ $\infty$, $\hcap(K(t))=2t$, and $\vphi(t,\cdot)=\vphi_{K(t)}$.
We call $K(t)$ and $\vphi(t,\cdot)$, $0\le t<T$, the chordal Loewner hulls and maps,
respectively, driven by $\xi$.

Let $B(t)$, $0\le t<\infty$, be a (standard) Brownian motion.
Let $\kappa > 0$. Then $K(t)$ and $\vphi(t,\cdot)$, $0\le
t<\infty$, driven by $\xi(t)=\sqrt\kappa B(t)$, $0\le t<\infty$, are
called the standard chordal SLE$(\kappa)$ hulls and maps, respectively.
It is known (\cite{RS-basic}\cite{LSW-2}) that almost surely for any
$t\in[0,\infty)$, \BGE\beta(t):=\lim_{\HH\ni
z\to\xi(t)}\vphi(t,\cdot)^{-1}(z)\label{trace}\EDE exists, and
$\beta(t)$, $0\le t<\infty$, is a continuous curve in $\lin{\HH}$.
Moreover, if $\kappa\in(0,4]$ then $\beta$ is a simple curve, which
intersects $\R$ only at the initial point, and for any $t\ge 0$,
$K(t)=\beta((0,t])$; if $\kappa>4$ then $\beta$ is not simple, and
intersects $\R$ at infinitely many points; and in general, $\HH\sem
K(t)$ is the unbounded component of $\HH\sem\beta((0,t])$ for any
$t\ge 0$. Such $\beta$ is called a standard chordal SLE$(\kappa)$
trace.

If $(\xi(t))$ is a semi-martingale with $d\langle \xi\rangle_t
=\kappa dt$ for some $\kappa>0$, then from the Girsanov's theorem (\cite{RY}) and the
existence of standard chordal SLE$(\kappa)$ trace, we see that almost surely for
any $t\in[0,T)$, $\beta(t)$ defined by (\ref{trace}) exists, and has
the same property as a standard chordal SLE$(\kappa)$ trace
(depending on the value of $\kappa$) as described in the last
paragraph.

Let $\kappa> 0$, $N\in\N$, $\vec{\rho}=(\rho_1,\dots,\rho_N)\in\R^N$, $x_0\in\R$, and
$\vec{p}=(p_1,\dots,p_N)\in(\ha\R\sem\{x_0\})^N$, where $\ha\R=\R\cup\{\infty\}$ is
a circle. Let $B(t)$ be a Brownian motion, which generates a filtration $(\F_t)$.
Let $\xi(t)$ and $p_m(t)$, $1\le m\le N$, $0\le t<T$, be the maximal solutions
to the SDE:
\begin{equation}\left\{\begin{array}{lll} d\xi(t) & = &
\sqrt\kappa d B(t)+\sum_{m=1}^N\frac{\rho_m\, dt}{\xi(t)-p_m(t)}\\ \\
dp_m(t) & = & \frac{2dt}{p_m(t)-\xi(t)},\quad 1\le m\le
N,\end{array}\right.\label{kappa-rho}\end{equation}
 with initial
values
$$ \xi(0)=x_0, \quad p_m(0)=p_m, \quad 1\le m\le N.$$
The meaning of the maximal solutions is that $[0,T)$ is the maximal interval of the solution.
Here if some $p_m=\infty$ then
$p_m(t)=\infty$ and $\frac{\rho_m}{\xi(t)-p_m(t)}=0$ for all $t\ge
0$, so $p_m$ has no effect on the equation.  Let $K(t)$, $0\le t<T$, be the chordal
Loewner hulls driven by $\xi$. Then we call $K(t)$, $0\le t<T$, a
 chordal SLE$(\kappa;\rho_1,\dots,\rho_N)$ or SLE$(\kappa;\vec{\rho})$  process started
from $(x_0;p_1,\dots,p_N)$ or $(x_0;\vec{p})$. It is known that $(\xi(t))$ is an $(\F_t)$-semi-martingale with
$d\langle \xi\rangle_t =\kappa dt$. So the chordal Loewner trace
$\beta(t)$, $0\le t<T$, driven by $\xi$ exists, and is called a chordal
 SLE$(\kappa;\vec{\rho})$ trace started from
 $(x_0;\vec{p})$. These $p_m$'s and $\rho_m$'s are called the force points and forces, respectively.

The chordal SLE$(\kappa;\vec{\rho})$ processes defined above are of generic
cases. We now introduce degenerate SLE$(\kappa;\vec{\rho})$ processes, where
one of the force points takes value $x_0^+$ or $x_0^-$, or two of the
force points take values $x_0^+$ and $x_0^-$, respectively. The force point $x_0^+$
or $x_0^-$ is called a degenerate force point. The definitions are as follows.
Suppose $p_1=x_0^+$ is the only degenerate force point. Let $\xi(t)$ and $p_m(t)$, $1\le
k\le N$, $0<t<T$, be the maximal solution to (\ref{kappa-rho}) with
initial values
$$ \xi(0)=p_1(0)=x_0, \quad p_k(0)=p_k, \quad 2\le k\le N.$$
Moreover, we require that
\BGE p_1(t)>\xi(t),\quad 0<t<T.\label{p>xi}\EDE
It is known that the solution exists, and $(\xi(t))$ is also an  $(\F_t)$-semi-martingale with
$d\langle \xi\rangle_t =\kappa dt$.
The chordal Loewner trace driven by $\xi(t)$, $0\le t<T$, is called a chordal
SLE$(\kappa;\rho_1,\dots,\rho_N)$ trace started from
$(x_0;x_0^+,p_2,\dots,p_N)$. If the ``$>$' in (\ref{p>xi}) is replaced by ``$<$'',
then we get a chordal SLE$(\kappa;\rho_1,\dots,\rho_N)$  trace started from
$(x_0;x_0^-,p_2,\dots,p_N)$. If the only degenerate force points are $p_1=x_0^+$ and $p_2=x_0^-$, let $\xi(t)$ and $p_k(t)$,
$1\le k\le N$, $0<t<T$, be the maximal solution to (\ref{kappa-rho})
with initial values $$ \xi(0)=p_1(0)=p_2(0)=x_0, \quad p_k(0)=p_k, \quad 3\le k\le N$$
such that  $$p_1(t)>\xi(t)>p_2(t),\quad 0<t<T.$$
The chordal Loewner trace driven by $\xi(t)$, $0\le t<T$, is called a
chordal SLE$(\kappa;\rho_1,\dots,\rho_N)$
trace started from $(x_0;x_0^+,x_0^-,p_3,\dots,p_N)$.

For $1\le m\le N$, the function $p_m(t)$, $0\le t<T$, is called the force point function started from $p_m$.
Each force point function is determined by its initial point $p_m$ and the driving function $\xi(t)$ as follows.
Let $\vphi(t,\cdot)$, $0\le t<T$, be the chordal Loewner maps driven by $\xi$. If $p_m$ is not degenerate,
then from (\ref{chordal-equation}), we have $p_m(t)=\vphi(t,p_m)$, $0\le t<T$. If $p_m=x_0^\sigma$, $\sigma
\in\{+,-\}$, is degenerate, then it is not difficult to see that $p_m(t)=\lim_{x\to x_0^\sigma}\vphi(t,x)$.

\vskip 3mm

The following lemma is a special case of Lemma 2.1 in \cite{duality}.

\begin{Lemma} Suppose $\kappa\in(0,4]$ and
$\vec{\rho}=(\rho_1,\dots,\rho_N)$ with  $\sum_{m=1}^N\rho_m =\kappa-6$.
For $j=1,2$, let $K_j(t)$, $0\le
t<T_j$, be a generic or degenerate chordal SLE$(\kappa;\vec{\rho})$ process started from
$(x_{j};\vec{p}_{j})$, where $\vec{p}_{j}=(p_{j,1},\dots,p_{j,N})$,
$j=1,2$. Suppose $W$ is a conformal
or conjugate conformal map from $\HH$ onto $\HH$ such that $W(x_1)=x_2$ and
$W(p_{1,m})=p_{2,m}$, $1\le m\le N$. Then $(W(K_1(t)),0\le t<T_1)$ has the same law as
$(K_2(t),0\le t<T_2)$ up to a time-change. A similar result holds
for the traces. \label{coordinate}
\end{Lemma}

The following lemma is a special case of Theorem 3.2 in \cite{duality}.

\begin{Lemma} Suppose $\kappa\in(0,4]$, $\rho\ge(\kappa-4)/2$, and
 $\beta(t)$, $0\le t<\infty$, is a chordal SLE$(\kappa; \rho)$
 trace started from $(0; 0^\sigma)$, where $\sigma\in\{+,-\}$. Then a.s.\ $\lim_{t\to\infty}\beta(t)=\infty$. \label{tendtoinfty}
\end{Lemma}

From Lemma \ref{coordinate} and Lemma \ref{tendtoinfty}, we obtain the following lemma.

\begin{Lemma} Let $\kappa\in(0,4]$, $\rho\ge (\kappa-4)/2$, and
 $x_1\ne x_2\in\R$. Suppose $\beta(t)$, $0\le t<T$, is a chordal SLE$(\kappa;\rho,\kappa-6-\rho)$ trace
 started from $(x_1;x_1^\sigma,x_2)$, where $\sigma\in\{+,-\}$. Then a.s.\ $\lim_{t\to T^-} \beta(t)=x_2$.
\label{??}
\end{Lemma}
{\bf Proof.} Let
 $\beta_0(t)$, $0\le t<\infty$, be a chordal SLE$(\kappa;\rho)$ trace started from
$(0;0^+)$. From Lemma \ref{tendtoinfty}, a.s.\ $\lim_{t\to\infty}\beta_0(t)=\infty$. We may find $W$
that maps $\HH$ conformally or conjugate conformally onto $\HH$ such that $W(0)=x_1$, $W(\infty)=x_2$,
and $W(0^+)=x_1^\sigma$. From Lemma \ref{coordinate}, after a time-change, $W(\beta_0(t))$, $0\le t<\infty$, has the same distribution
as $\beta(t)$, $0\le t<T$. Thus, a.s.\ $\lim_{t\to T^-}\beta(t)=W(\infty)=x_2$.  $\Box$

\section{Intermediate SLE$(\kappa;\rho)$  Process}

\begin{Lemma} For $\kappa\in(0,4)$ and $\rho\ge (\kappa-4)/2$, let
$a=\frac{2\rho}\kappa$, $b=1-\frac4\kappa<0$, and $c= \frac{2\rho+4}\kappa\ge 1$.  For
$x\in(-1,1)$, let $U_0(x)={_2F_1}(a,b;c;x)$, where $_2F_1$ is the hypergeometric function \cite{hand}. Then
there are $C_2>C_1>0$ such that $C_1\le U_0(x)\le C_2$ on $[0,1)$. Let $f_0(x)=\frac{U_0'(x)}{U_0(x)}$ on $[0,1)$. Then $f_0$ is also bounded
on $[0,1)$, $f_0(x)\ge \frac{b}{1-x}$ for $0\le x<1$, and $\lim_{x\to 1^-}f_0(x)=-\frac a2$.\label{f>*}\end{Lemma}
{\bf Proof.} It is known \cite{hand} that $U_0$ is analytic and satisfies the Gaussian hypergeometric equation:
\BGE x(x-1)U_0''(x)+[(a+b+1)x-c]U_0'(x)+ab U_0(x)=0.\label{Gauss-hyper*}\EDE
Moreover, we have $U_0(0)=1> 0$ and $f_0'(0)=U_0'(0)=\frac{ab}c$.
Let $z_0=\sup\{x\in(0,1):U_0(x)\ne0\}$. Then $z_0\in(0,1]$ and $f_0$ is analytic on $[0,z_0)$.
Let $h_0(x)=f_0(x)-\frac{b}{1-x}=\frac{U_0'(x)}{U_0(x)}-\frac{b}{1-x}$ on $[0,z_0)$. Then
$h_0(0)=\frac{ab}c -b=\frac{-4b}{2\rho+4}>0$.
From (\ref{Gauss-hyper*}) and that $b+c-a=1$, we find that for $x\in[0,z_0)$, $h_0(x)$ satisfies
\BGE x h_0'(x)+x h_0(x)^2+ch_0(x)+\frac{b(1-b)}{(1-x)^2}=0.\label{h}\EDE
Assume that there is $x_1\in[0,z_0)$ such that $h_0(x_1)\le 0$. Since $h_0(0)>0$, so $x_1>0$ and there is $x_0\in(0,z_0)$ such that
$h_0(x_0)=0$ and $h_0(x)>0$ for $x\in[0,x_0)$. Then we have $h_0'(x_0)\le 0$. However, since $b<0$, from (\ref{h}) we have $h_0'(x_0)>0$, which
is a contradiction. Thus $h_0(x)>0$ for all $x\in[0,z_0)$. So we have $f_0(x)> \frac{b}{1-x}$ for $0\le x<z_0$.
Assume that $z_0<1$. Then $z_0$ is a zero of $U_0$, so $z_0$ is a simple pole of
$f_0$, and the residue is positive. Thus, $\lim_{x\to z_0^-}f_0(x)=-\infty$, which contradicts that
$f_0(x)> \frac{b}{1-x}$ for $0\le x<z_0$. Thus, $z_0=1$. So $U_0(x)\ne 0$ and
$f_0(x)> \frac{b}{1-x}$ for $0\le x<1$. Since $U_0(0)=1>0$, so $U_0(x)>0$ on $[0,1)$.

Now $U_0$ and $f_0$ are continuous on $[0,1)$, and $U_0(x)>0$ on $[0,1)$. To complete the proof, we suffice to show
that  $\lim_{x\to 1^-}U_0(x)$ and $\lim_{x\to 1^-}f_0(x)$ both exist and are finite, and $\lim_{x\to 1^-}U_0(x)> 0$.
One may check that $c$, $c-a$, $c-b$ and $c-a-b$ are all positive. So from \cite{hand},
\BGE \lim_{x\to 1^-}U_0(x)=\frac{\Gamma(c)\Gamma(c-a-b)}{\Gamma(c-a)\Gamma(c-b)}\in(0,\infty). \label{u1*}\EDE
We have $U_0'(x)=\frac{ab}c  {_2F_1}(a+1,b+1;c+1;x)$. One may check that $c+1$ and $(c+1)-(a+1)-(b+1)$ are both positive.
So from \cite{hand} again,
\BGE \lim_{x\to 1^-} U_0'(x)=\frac{ab}c \cdot \frac{\Gamma(c+1)\Gamma(c-a-b-1)}{\Gamma(c-a)\Gamma(c-b)}.\label{u'1*}\EDE
From (\ref{u1*}) and (\ref{u'1*}), we have $\lim_{x\to 1^-} f_0(x)=\frac{ab}{c-a-b-1}=-\frac{a}2$, which is finite. $\Box$

\vskip 3mm

From now on, fix $\kappa\in(0,4)$ and $\rho\ge(\kappa-4)/2$. Let $f_0$ be given
by Lemma \ref{f>*}. Let
\BGE g_0(x):=\rho+\kappa x f_0(x).\label{g0}\EDE
From Lemma \ref{f>*}, $g_0$ is bounded on $[0,1)$, $\lim_{x\to 1^-}g_0(x)=0$, and for $0\le x<1$,
\BGE g_0(x)\ge \rho+(\kappa-4)\frac{x}{1-x}.\label{g0>}\EDE
For $0<p_1<p_2$, let
\BGE J(p_1,p_2):=-\Big(\frac 1{p_1}-\frac 1{p_2}\Big)g_0\Big(\frac{p_1}{p_2}\Big).\label{J}\EDE
From (\ref{g0>}) and that $\rho\ge \kappa/2-2$, we have
\BGE J(p_1,p_2)\le \frac{\rho}{p_2}-\frac{\rho}{p_1}+\frac{4-\kappa}{p_2}\le
\frac{2-\kappa/2}{p_1}+\frac{2-\kappa/2}{p_2}.\label{R<*}\EDE

Let $0<p_1<p_2$. Let $B(t)$ be a Brownian motion. Let $J(\cdot,\cdot)$ be defined by (\ref{J}).
Let $\xi(t)$, $p_1(t)$ and $p_2(t)$, $0\le t<T$, be the maximal solution to
\BGE
\left\{ \begin{array}{l}
d\xi(t)=\sqrt\kappa dB(t)+J(p_1(t)-\xi(t),p_2(t)-\xi(t))dt, \\ \\
dp_1(t)=\frac {2dt}{p_1(t)-\xi(t)},\quad dp_2(t)=\frac {2dt}{p_2(t)-\xi(t)},
\end{array} \right.\label{generic-eqn*}
\EDE
with initial values
$$\xi(0)=0, \qquad p_j(0)=p_j, \quad j=1,2.$$
We call
the chordal Loewner trace $\beta(t)$, $0\le t<T$, driven by $\xi$,  a (generic) intermediate
SLE$(\kappa;\rho)$ trace with force points $p_1$ and $p_2$.
Note that $\xi(t)<p_1(t)<p_2(t)$ for $0\le t<T$. If $T<\infty$, we must have $\lim_{t\to T^-}
p_1(t)-\xi(t)=0$. Thus, if $\limsup_{t\to T^-}p_1(t)-\xi(t)>0$, then $T=\infty$.

\begin{Theorem}  Let $\beta(t)$, $0\le t<T$, be an intermediate SLE$(\kappa;\rho)$ trace. Then
a.s.\ $T=\infty$,
which means that $\infty$ is a subsequential limit of $\beta(t)$ as $t\to T^-$. \label{T=infty*}
\end{Theorem}
{\bf Proof.} Let $\xi(t)$, $0\le t<T$, be the driving function for $\beta$. Then there are $p_1(t)$, $p_2(t)$
and some Brownian motion $B(t)$ such that (\ref{generic-eqn*}) holds, and $[0,T)$ is the maximal interval of
the solution. Let $X_j(t)=p_j(t)-\xi(t)$, $j=1,2$. Then $0<X_1(t)<X_2(t)$, $0\le t<T$; and
for $j=1,2$, $X_j$ satisfies the SDE
$$dX_j(t)=-\sqrt\kappa dB(t)+\Big(\frac{2}{X_j(t)}-J(X_1(t),X_2(t))\Big)\,dt.$$
From It\^o's formula (\cite{RY}), for $j=1,2$, we have
\BGE d\ln(X_j(t))=-\frac{\sqrt\kappa}{X_j(t)}\,dB(t)+\Big(\frac{2-\kappa/2}{X_j(t)^2}-\frac{J(X_1(t),X_2(t))}{X_j(t)}\Big)\,dt.\label{Xj*}\EDE
Thus, we have
$$d(\ln(X_2(t)/X_1(t)))=\Big(\frac{\sqrt\kappa}{X_1(t)}-\frac{\sqrt\kappa}{X_2(t)}\Big)\,dB(t)
-\Big(\frac{2-\kappa/2}{X_1(t)^2}-\frac{2-\kappa/2}{X_2(t)^2}\Big)\,dt$$
$$+\Big(\frac{1}{X_1(t)}-\frac{1}{X_2(t)}\Big)J(X_1(t),X_2(t))\,dt.$$
Since $1/X_1(t)>1/X_2(t)$ and $2-\kappa/2>0$, so
from (\ref{R<*}),  the drift term for $\ln(X_2(t)/X_1(t))$ is not positive. Note that $\ln(X_2(t)/X_1(t))$ is always positive.
So $(\ln(X_2(t)/X_1(t)))$  is a supermartingale. Thus, a.s.\
$\lim_{t\to T^-}\ln(X_2(t)/X_1(t))$ exists and is finite. So a.s.\
\BGE\int_0^T \Big(\frac{\sqrt\kappa}{X_1(t)}-\frac{\sqrt\kappa}{X_2(t)}\Big)^2\,dt=\lim_{t\to T^-}\langle \ln(X_2/X_1)\rangle_t<\infty.
\label{integ<infty1*}\EDE

Let ${\cal E}_1$ denote the event that $\lim_{t\to T^-}\ln(X_2(t)/X_1(t))> 0$. Assume that ${\cal E}_1$ occurs.
 From (\ref{integ<infty1*}),
we have a.s.\ $\int_0^T X_1(t)^{-2}dt<\infty$. From (\ref{J}) and (\ref{Xj*}), we have
\BGE d\ln(X_1(t))=-\frac{\sqrt\kappa}{X_1(t)}\,dB(t)+\frac{1}{X_1(t)^2}\Big[{2-\frac\kappa2}+
\Big(1-\frac{X_1(t)}{X_2(t)}\Big) \,g_0
\Big(\frac{X_1(t)}{X_2(t)}\Big)\Big]\,dt.\label{X1*}\EDE
Since a.s.\ $\int_0^T X_1(t)^{-2}dt<\infty$, and $g_0$ is bounded on $[0,1)$, so a.s.\
$$\int_0^T \frac{1}{X_1(t)^2} \Big|{2-\frac\kappa2}+
\Big(1-\frac{X_1(t)}{X_2(t)}\Big) \,g_0
\Big(\frac{X_1(t)}{X_2(t)}\Big)\Big|\,dt<\infty.$$
From (\ref{X1*}) we have a.s.\ $\lim_{t\to T^-} \ln(X_1(t))$ exists and is finite. Thus, on ${\cal E}_1$
a.s.\ $\lim_{t\to T^-}  X_1(t)$ exists and
is positive, which implies that $T=\infty$.

Let ${\cal E}_2$ denote the event that $\lim_{t\to T^-}\ln(X_2(t)/X_1(t))= 0$. Assume that ${\cal E}_1$ occurs.
Then $\lim_{t\to T^-} X_1(t)/X_2(t)=1$, so
$\lim_{t\to T^-} g_0(X_1(t)/X_2(t))=\lim_{x\to 1^-}g_0(x)=0$. Since $2-\kappa/2>0$,
so the drift term in (\ref{X1*}) is
positive when $t$ is close to $T$. From (\ref{X1*}), a.s.\ $\limsup_{t\to T^-} \ln(X_1(t))>-\infty$, which implies that
$\limsup_{t\to T^-}  X_1(t)>0$. So we have a.s.\ $T=\infty$ on the event ${\cal E}_2$.

Since ${\cal E}_1\cup{\cal E}_2$ is a.s.\ the whole probability space, so a.s.\ $T=\infty$. Suppose $T=\infty$.
Since for any $0<t<\infty$, the half-plane capacity of $\beta((0,t])$ is $2t$, so the diameter of $\beta((0,t])$ is at least $\sqrt{2t}$.
Thus, the diameter of $\beta((0,\infty))$ is infinite, so $\infty$ is a subsequential limit of $\beta(t)$ as $t\to T^-$. $\Box$

\vskip 3mm

The above theorem still holds if the force points $p_1$ and $p_2$ are random points, and the joint distribution of $p_1$ and $p_2$
is independent of the Brownian motion $B(t)$. The argument in the above proof still works.

\vskip 3mm

We may let the force point $p_1$ be $0^+$, and define the degenerate intermediate SLE$(\kappa;\rho)$ trace.
The definition is as follows. Fix $p_2>0$. Let $\xi(t)$, $p_1(t)$ and $p_2(t)$ solve (\ref{generic-eqn*}) for $0<t<T$,
with initial values \BGE \xi(0)=p_1(0)=0,\qquad p_2(0)=p_2.\label{initial-int}\EDE Moreover, we require that
\BGE\xi(t)<p_1(t), \qquad 0<t<T.\label{ineq-int}\EDE
The chordal Loewner trace $\beta(t)$, $0\le t<T$, driven by $\xi$, is called a degenerate intermediate
SLE$(\kappa;\rho)$ trace  with force points $0^+$ and $p_2$.

We claim that the solution to (\ref{generic-eqn*}) together with (\ref{initial-int}) and (\ref{ineq-int})
a.s.\ exists. For the proof, we suffice to prove that the
solution exists on $(0,T_0)$ for some stopping time $T_0>0$ because after $T_0$ we are dealing with some generic
case with random force points. Let
$\til B(t)$ be a Brownian motion under some probability measure $\PP$. Let $\xi(t)$, $p_1(t)$ and $p_2(t)$,
$0< t<T_1$, be the maximal solution to $$
\left\{ \begin{array}{l}
d\xi(t)=\sqrt\kappa d\til B(t)+\frac{\rho }{\xi(t)-p_1(t)}\,dt, \\
dp_j(t)=\frac {2dt}{p_j(t)-\xi(t)},\quad j=1,2,
\end{array} \right.
$$
such that (\ref{initial-int}) and (\ref{ineq-int}) hold.
The solution a.s.\ exists because $\xi$ is the driving function
for an SLE$(\kappa;\rho)$ process started from $(0,0^+)$.

From (\ref{g0}) and (\ref{J}), it is clear that $\lim_{p_1\to 0^+} \Big(J(p_1,p_2)+\frac{\rho}{p_1}\Big)
=\frac\rho{p_2}-\frac\kappa{p_2}\,f_0(0)$.
Define $Z(t)$, $0\le t<T_1$, such that for $t>0$, $Z(t)=J(p_1(t)-\xi(t),p_2(t)-\xi(t))-\frac{\rho}{\xi(t)-p_1(t)}$, and
$Z(0)=\frac\rho{p_2}-\frac\kappa{p_2}\,f_0(0)$. Then
$Z(t)$ is continuous on $[0,T_1)$. From the Girsanov's Theorem, there is a stopping
time $T_0\in(0,T_1)$ such that under some other probability measure $\QQ$, $B(t):=\til B(t)-\frac1{\sqrt\kappa}\int_0^t Z(s)ds$, $0\le t<T_0$, is
a partial Brownian motion, which means that $B(t)$ could be extended to a full Brownian motion. Then we have
$$d\xi(t)=\sqrt\kappa dB(t)+J(p_1(t)-\xi(t),p_2(t)-\xi(t))dt,\quad 0\le t<T_0.$$
Thus, the solution to (\ref{generic-eqn*}) with (\ref{initial-int}) and (\ref{ineq-int}) a.s.\
exists on $(0,T_0)$. Then the solution can be
extended to the maximal interval, say $(0,T)$, and so we have the existence of the maximal solution.
From Theorem \ref{T=infty*}, we get the following  corollary.

\begin{Corollary} Let $\beta(t)$, $0\le t<T$, be a degenerate intermediate SLE$(\kappa;\rho)$ trace.
Then a.s.\ $T=\infty$, which means that $\infty$ is a subsequential limit of $\beta(t)$ as $t\to T^-$.
\label{T=infty-degenerate*}
\end{Corollary}

\section{Martingales} \label{martingale}
Fix $\kappa\in(0,4)$ and $\rho\ge \kappa/2-2$. Let $x_1<x_2\in\R$, $\sigma_1=+$ and $\sigma_2=-$.
Throughout this section, the subscripts $j$ and $k$ will be any of the two numbers: $1$ or $2$, such
that $j$ and $k$ are different. Let $\xi_j(t)$, $0\le t<T_j$, be the driving function for a
chordal SLE$(\kappa;\rho,\kappa-6-\rho)$ trace $\beta_j(t)$, $0\le t<T_j$,
 started from $(x_j;x_j^{\sigma_j},x_{k})$.
 From Lemma \ref{??}, we have a.s.\ $\lim_{t\to T_j^-} \beta_j(t)=x_k$.
Let $\vphi_j(t,\cdot)$ and $K_j(t)$, $0\le t<T_j$, be the chordal Loewner maps and hulls driven by $\xi_j$.
 Let $p_j(t)$ and $q_j(t)$ be the force point functions started from $x_j^{\sigma_j}$ and $x_k$, respectively.
So we have $p_j(t)=\lim_{x\to x_j^{\sigma_j}}\vphi_j(t,x)$ and $q_j(t)=\vphi_j(t,x_k)$.
For $0\le t<T$, let
$$B_j(t)=\frac 1{\sqrt\kappa}\,\Big(\xi_j(t)-x_j-\int_0^t\frac{\rho}{\xi_j(s)-p_{j}(s)}\,ds+\int_0^t
\frac{\kappa-6-\rho}{\xi_j(s)-q_{j}(s)}\,ds\Big).$$
Then $B_j(t)$, $0\le t<T$, is a partial Brownian motions. Let $(\F^j_t)$ be the filtration generated by $B_j(t)$.
Then $(\xi_j(t))$, $p_j(t)$, and $(q_j(t))$ are all $(\F_t)$-adapted. And $(\xi_j(t))$ is an  $(\F_t)$-semi-martingale
with $d\langle \xi\rangle_t=\kappa dt$. Moreover, $\xi_j(t)$, $p_{j}(t)$ and $q_{j}(t)$, $0<t<T_j$, are
the maximal solution to the following equations
\begin{eqnarray}
d\xi_j(t)& = &\sqrt\kappa dB_j(t)+\frac{\rho}{\xi_j(t)-p_{j}(t)}\,dt+\frac{\kappa-6-\rho}{\xi_j(t)-q_{j}(t)}\,dt,\label{kappa-rho-deg*-1}\\
dp_{j}(t)&=&\frac{2}{p_{j}(t)-\xi_j(t)}\,dt,\label{kappa-rho-deg*-2}
\\  dq_{j}(t)&=&\frac{2}{q_{j}(t)-\xi_j(t)}\,dt\label{kappa-rho-deg*-3},
\end{eqnarray}
with initial values
\BGE \xi_j(0)=p_{j}(0)=x_j, \quad q_{j}(0)=x_{k};\label{initial} \EDE
and they satisfy the inequalities
\BGE \xi_1(t)<p_1(t)<q_1(t),\quad 0<t<T_1;\qquad \xi_2(t)>p_2(t)>q_2(t),\quad 0<t<T_2.\label{ineq} \EDE

Now suppose that $(\xi_1(t))$ and $(\xi_2(t))$ are independent. Then $(B_1(t))$ and $(B_2(t))$ are also independent.
So for any fixed $(\F^k_t)$-stopping time $t_k$ with $0\le t_k<T_k$, $B_j(t)$, $0\le t<T_j$, is a partial
$(\F^j_t\times\F^k_{t_k})_{t\ge 0}$-Brownian motion.

Differentiating (\ref{chordal-equation}) w.r.t.\ $\pa_z$ and plugging $\xi=\xi_j$ and $z=x_k$,
we find that for $0\le t<T_j$,
\BGE \frac{d\pa_z \vphi_j(t,x_k)}{\pa_z \vphi_j(t,x_k)}=\frac{-2 dt}{(q_j(t_j)-\xi_j(t_j))^2}.\label{q2*}\EDE
From (\ref{kappa-rho-deg*-1})-(\ref{kappa-rho-deg*-3}) we have that, for $0<t<T_j$,
\begin{eqnarray}
\frac{d
(\xi_j(t)-p_j(t))}{\xi_j(t)-p_j(t)}&=&\frac{ d\xi_j(t)}{\xi_j(t)-p_j(t)}+\frac{2d t}{(\xi_j(t)-p_j(t))^2};
\label{jjj*}\\
\frac{d
(\xi_j(t)-q_j(t))}{\xi_j(t)-q_j(t)}&=&\frac{ d\xi_j(t)}{\xi_j(t)-q_j(t)}+ \frac{2 d t}{(\xi_j(t)-q_j(t))^2};
\label{jjk*}\\
\frac{d
(q_j(t)-p_j(t))}{q_j(t)-p_j(t)}&=&\frac{-2 d t}{(\xi_j(t)-q_j(t))(\xi_j(t)-p_j(t))}.\label{jkj*}
\end{eqnarray}
In the above equations, (\ref{q2*}) and (\ref{jkj*}) are ODEs, (\ref{jjj*}) and (\ref{jjk*}) are $(\F^j_t)$-adapted SDEs.

For    $t\in(0,T_j)$,
define \BGE r_j(t)=|\xi_j(t)-p_j(t)|^{-\frac\rho\kappa}|\xi_j(t)-q_j(t)|^{-\frac{\kappa-6-\rho}{\kappa}}
|q_j(t)-p_j(t)|^{-\frac{\rho(\kappa-6-\rho)}{2\kappa}}\pa_z \vphi_j(t,x_{k})^{\frac{(\rho+2)(\kappa-6-\rho)}{4\kappa}}.\label{r*}\EDE
From (\ref{kappa-rho-deg*-1}), (\ref{q2*})-(\ref{jkj*}) and It\^o's formula, we have that, for $t>0$,
\BGE \frac{dr_j(t)}{r_j(t)}=-\frac{\rho}{\xi_j(t)-p_j(t)}\cdot\frac{dB_j(t)}{\sqrt\kappa} -\frac{{\kappa-6-\rho} }{\xi_j(t)-q_j(t)}\cdot
\frac{dB_j(t)}{\sqrt\kappa} +\frac{{\rho(\kappa-4-\rho)}/({2\kappa})}{(\xi_j(t)-p_j(t))^2}\,d t.\label{ODE-r*}\EDE

Let ${\cal D}=\{(t_1,t_2)\in[0,T_1)\times[0,T_2):\beta_1([0,t_1])\cap\beta_2([0,t_2])=\emptyset\}$.
Then for any $(t_1,t_2)\in\cal D$, $K_1(t_1)\cup K_2(t_2)$ is a hull in $\HH$ w.r.t.\ $\infty$.
For $(t_1,t_2)\in{\cal D}$, let \BGE
K_{k,t_j}(t_k):=(K_j(t_j)\cup
K_{k}(t_k))/K_j(t_j)=\vphi_{j}(t_j,K_k(t_k)),\label{K-j-t*}\EDE and $\vphi_{k,t_j}(t_k,\cdot):=\vphi_{K_{k,t_j}(t_k)}$.
Then $K_{k,t_j}(t_k)$ is the image of a curve in $\HH$ started from $\vphi_j(t_j,x_k)=q_j(t_j)$.
And for any $z\in\HH\sem(K_1(t_1)\cup K_2(t_2))$, \BGE\vphi_{K_1(t_1)\cup
K_2(t_2)}(z)=\vphi_{1,t_2}(t_1,\vphi_2(t_2,z)) =
\vphi_{2,t_1}(t_2,\vphi_1(t_1,z)) .\label{circ=1*}\EDE

Define $A_{j,h}$, $h\in\Z_{\ge 0}$, on $\cal D$ such that
$A_{j,h}(t_1,t_2)=\pa_z^h\vphi_{k,t_j}(t_k,\xi_j(t_j))$. Note that the definition of $A_{j,h}$ here
agrees with the definition of $A_{j,h}$ in Section 4.2 of \cite{duality}.
From now on, we fix $t_k$ to be some $(\F^k_t)$-stopping time that lies on $[0,T_k)$, and
consider the filtration $(\F^j_{t_j}\times\F^k_{t_k})_{t_j\ge 0}$. Since $B_j(t)$ and $B_k(t)$ are independent
Brownian motions, so $B_j(t_j)$ is an $(\F^j_{t_j}\times\F^k_{t_k})_{t_j\ge 0}$-Brownian motion.
We use $\pa_j$ to denote the partial derivative w.r.t.\  $t_j$.
The following  equations are (4.10) and (4.12) in \cite{duality}, where
(\ref{pa-A_1-0}) is an $(\F^j_{t_j}\times\F^k_{t_k})_{t_j\ge 0}$-adapted SDE.
\BGE \pa_j A_{j,0}= A_{j,1} \pa \xi_j(t_j)+(\frac{\kappa} 2-3)A_{j,2}\pa t_j;\label{pa-A_1-0}\EDE
\BGE \pa_j A_{k,0}=\frac{2
A_{j,1}^2}{A_{k,0}-A_{j,0}},\qquad \frac{\pa_j
A_{k,1}}{A_{k,1}}=\frac{-2 A_{j,1}^2}{(A_{k,0}-A_{j,0})^2}.\label{pa-A_2}\EDE

We now use $\pa_1$ and $\pa_z$ to denote the partial
derivatives of
$\vphi_{j,t_0}(\cdot,\cdot)$ w.r.t.\ the first (real) and second
(complex) variables, respectively, inside the bracket; and use
$\pa_0$ to denote the partial derivative of
$\vphi_{j,t_0}(\cdot,\cdot)$ w.r.t. the subscript $t_0$. Let $(t_1,t_2)\in\cal D$.
The following equations are (3.9) and (3.15) in \cite{reversibility}.
\begin{eqnarray}\pa_1\vphi_{j,t_k}(t_j,z)&=&\frac{2A_{j,1}^2}
{\vphi_{j,t_k}(t_j,z)-A_{j,0}},\quad z\in \HH\sem K_{j,t_k}(t_j);\label{chordal*-}\\
  \pa_0\vphi_{k,t_j}(t_k,z)&=&\frac{2A_{j,1}^2}
{\vphi_{k,t_j}(t_k,z)-A_{j,0}}-\frac{2\pa_z\vphi_{k,t_j}(t_k,z)}
{z-\xi_j(t_j)},\quad z\in \HH\sem K_{k,t_j}(t_k).\label{pa0-1*}\end{eqnarray}
Since $\lin{K_{j,t_k}(t_j)}\cap\R=\{q_k(t_k)\}$ and $\lin{ K_{k,t_j}(t_k)}\cap\R=\{q_j(t_j)\}$, so
after continuation, (\ref{chordal*-}) also holds for any $z\in \R\sem \{q_k(t_k)\}$, and
(\ref{pa0-1*}) also holds for any $z\in \R\sem \{\xi_j(t_j),q_j(t_j)\}$.
Differentiating (\ref{pa0-1*}) w.r.t.\ $\pa_z$, we find that for  $(t_1,t_2)\in\cal D$,
and $z\in \R\sem \{\xi_j(t_j),q_j(t_j)\}$,
\BGE \pa_0\pa_z\vphi_{k,t_j}(t_k,z)=-\frac{2A_{j,1}^2\pa_z\vphi_{k,t_j}(t_k,z)}
{(\vphi_{k,t_j}(t_k,z)-A_{j,0})^2}-\frac{2\pa_z^2\vphi_{k,t_j}(t_k,z)}{z-\xi_j(t_j)}
+\frac{2\pa_z\vphi_{k,t_j}(t_k,z)} {(z-\xi_j(t_j))^2}.\label{pa0-2*}\EDE

Define $B_{j,0}$ on $\cal D$ such that
$B_{j,0}(t_1,t_2)=\vphi_{k,t_j}(t_k,p_j(t_j))$.
Since $\xi_1(0)=p_1(0)$ and $\xi_1(t)<p_1(t)$ for $t>0$, so
$A_{1,0}(0,t_2)=B_{1,0}(0,t_2)$ and $A_{1,0}(t_1,t_2)<B_{1,0}(t_1,t_2)$ if $t_1>0$. Similarly, we have
$A_{2,0}(t_1,0)=B_{2,0}(t_1,0)$ and $A_{2,0}(t_1,t_2)>B_{2,0}(t_1,t_2)$ if $t_2>0$.
Choose any $y_1<y_2\in(x_1,x_2)$. Then $p_1(t_1)\le \vphi_1(t_1,y_1)<\vphi_2(t_1,y_2)$ for any $t_1\in[0,T_1)$.
From (\ref{circ=1*}) we have $$B_{1,0}(t_1,t_2)\le \vphi_{K_1(t_1)\cup K_2(t_2)}(y_1)<\vphi_{K_1(t_1)\cup K_2(t_2)}(y_2).$$
for any $(t_1,t_2)\in\cal D$. Similarly,
$B_{2,0}(t_1,t_2)\ge \vphi_{K_1(t_1)\cup K_2(t_2)}(y_2)>\vphi_{K_1(t_1)\cup K_2(t_2)}(y_1)$
for any $(t_1,t_2)\in\cal D$. Thus, $B_{1,0}<B_{2,0}$ on $\cal D$. So in general, $A_{1,0}\le B_{1,0}<B_{2,0}\le A_{2,0}$, where
$A_{1,0}=B_{1,0}$ iff $t_1=0$, and $B_{2,0}=A_{2,0}$ iff $t_2=0$.

Let   $(t_1,t_2)\in\cal D$. Since $p_k(t_k)\ne q_k(t_k)$, so we may
apply (\ref{chordal*-}) with $z= p_k(t_k)$, and obtain
\BGE \pa_j B_{k,0}=\frac{2A_{j,1}^2}{B_{k,0}-A_{j,0}}.\label{Bj0-1*}\EDE
Now suppose $t_j>0$. Then $p_j(t_j)\in\R\sem\{\xi_j(t_j),q_j(t_j)\}$. So we may
apply (\ref{pa0-1*}) with $z=p_j(t_j)$, and use (\ref{kappa-rho-deg*-2}) and chain rule to obtain
\BGE \pa_j B_{j,0}=\frac{2A_{j,1}^2}
{B_{j,0}-A_{j,0}}.\label{Bj0-2*}\EDE
Note that (\ref{Bj0-1*}) and (\ref{Bj0-2*}) have the same forms as the formula for $\pa_j B_{m,0}$ in
 (4.13) in \cite{duality}. But here we require that $t_j>0$ in (\ref{Bj0-2*}).

Let $E_{j,0}=A_{j,0}-A_{k,0}=-E_{k,0}\ne 0$,
$E_{j,m}=A_{j,0}-B_{m,0}$, $m=1,2$, and $C_{j,k}=B_{j,0}-B_{k,0}=-C_{k,j}\ne 0$.
From (\ref{pa-A_1-0})-(\ref{pa-A_2}) and (\ref{Bj0-1*})-(\ref{Bj0-2*}), we obtain the following formulas,
which have the same forms as (4.14) and (4.15) in \cite{duality}.
\begin{eqnarray}
\frac{\pa_j
E_{j,m}}{E_{j,m}}&=&\frac{A_{j,1}}{E_{j,m}}\,\pa\xi_j(t_j)+\Big(
\Big(\frac{\kappa}2-3\Big)\cdot\frac{A_{j,2}}{E_{j,m}}+2\cdot\frac{A_{j,1}^2}{E_{j,m}^2}\Big)\,\pa
t_j,\quad m=0,1,2;\label{Djm*}\\ \frac{\pa_j
E_{k,m}}{E_{k,m}}&=&\frac{-2A_{j,1}^2}{E_{j,0}E_{j,m}}\,\pa t_j,\quad m=1,2;\label{Dkm*} \\
\frac{\pa_j
C_{j,k}}{C_{j,k}}&=&\frac{-2A_{j,1}^2}{E_{j,1}E_{j,2}}\,\pa
t_j.\label{Cjk*}\end{eqnarray}
Here we require that $t_j>0$ in the SDEs for $\pa_j E_{j,j}$, $\pa_j E_{k,j}$, and  $\pa_j C_{j,k}$,
 because (\ref{Bj0-2*}) does not hold for $t_j=0$.

Define $\til B_{j,1}$ on $\cal D$ such that $\til B_{j,1}(t_1,t_2)=\pa_z\vphi_{k,t_j}(t_k,p_j(t_j))$.
Differentiating (\ref{chordal*-}) w.r.t.\ $\pa_z$ and plugging $z=p_k(t_k)$, we get
\BGE \frac{\pa_j \til B_{k,1}}{\til B_{k,1}}=\frac{-2 A_{j,1}^2}{E_{j,k}^2}\,\pa t_j.\label{Bj1-1*}\EDE
Applying (\ref{pa0-2*}) with $z=p_j(t_j)$, and using (\ref{kappa-rho-deg*-2}) and chain rule, we find that, for $t_j>0$,
\BGE \frac{\pa_j \til B_{j,1}}{\til B_{j,1}}=\frac{-2 A_{j,1}^2}{E_{j,j}^2}\,\pa t_j+\frac{2}{(p_j(t_j)-\xi_j(t_j))^2}\,\pa t_j.\label{Bj1-2*}\EDE
Let $D=\frac{\til B_{1,1}\til B_{2,1}}{C_{1,2}^2}=\frac{\til B_{1,1}\til B_{2,1}}{C_{2,1}^2}$.
From (\ref{Cjk*})-(\ref{Bj1-2*}), we find that, for $t_j>0$,
\BGE \frac{\pa_j D}{D}=-2\Big(\frac{A_{1,1}}{E_{j,j}}-\frac{A_{1,1}}{E_{j,k}}\Big)^2\,\pa t_j+\frac{2}{(p_j(t_j)-\xi_j(t_j))^2}\,\pa t_j.
\label{ODE-D*}\EDE

Let ${\cal D}'=\{(t_1,t_2)\in{\cal D}:t_1*t_2\ne 0\}$.
Define $R$ on $\cal D$ such that $R=\frac{E_{1,1}E_{2,2}}{E_{1,2}E_{2,1}}=\frac{E_{j,j}E_{k,k}}{E_{j,k}E_{k,j}}$.
From $A_{1,0}\le B_{1,0}<B_{2,0}\le A_{2,0}$ we have $|E_{j,j}|<|E_{j,k}|$ and $E_{j,j}/E_{j,k}\ge 0$,
so $R\in[0,1)$.  Since $A_{j,0}\ne B_{j,0}$ when $t_j>0$, so
$E_{1,1}*E_{2,2}\ne 0$ on ${\cal D}'$. Thus, $R\in(0,1)$ on ${\cal D}'$.
Since $E_{k,m}=E_{j,m}-E_{j,0}$ for $m=1,2$, so we have
\BGE \frac{R+1}{R-1}=\frac{2/E_{j,0}}{1/E_{j,j}-1/E_{j,k}}-\frac{1/E_{j,j}+1/E_{j,k}}{1/E_{j,j}-1/E_{j,k}}.\label{1+R}\EDE
From (\ref{Djm*}) and (\ref{Dkm*}), we have that, for $t_j>0$,
\begin{eqnarray}{\pa_j R}&=&R\Big(\frac{A_{j,1}}{E_{j,j}}-\frac{A_{j,1}}{E_{j,k}}\Big)\,\pa\xi_j(t_j)+R
\Big[\Big(\frac\kappa 2-3\Big)\Big(\frac{A_{j,2}}{E_{j,j}}-\frac{A_{j,2}}{E_{j,k}}\Big)
+\frac\kappa 2\,\Big(\frac{A_{j,1}}{E_{j,j}}-\frac{A_{j,1}}{E_{j,k}}\Big)^2\nonumber\\
& +&\Big(2-\frac\kappa 2\Big) \Big(\frac{A_{j,1}^2}{E_{j,j}^2}-\frac{A_{j,1}^2}{E_{j,k}^2}\Big)
+\Big(\frac{2A_{j,1}^2}{E_{j,0} E_{j,j}}-\frac{2A_{j,1}^2}{E_{j,0} E_{j,k}}\Big)
\Big]\,\pa t_j.\label{SDE-X*}\end{eqnarray}
Let $U_0(x)$ and $f_0(x)$ be given by Lemma \ref{f>*}. Let $g_0$ be defined
by (\ref{g0}). For $x\in(0,1)$, let $V_0(x):=x^{\frac\rho\kappa} U_0(x)$.
From (\ref{Gauss-hyper*}) and (\ref{g0}), we find that $V_0(x)$ satisfies
 \BGE x\frac{V_0'(x)}{V_0(x)}=\frac{g_0(x)}\kappa. \label{V0'}\EDE
\BGE \frac{\kappa}2\,\frac{V_0''(x)}{V_0(x)}x^2=\Big[\Big(2-\frac\kappa 2\Big)\frac{x+1}{x-1}-\frac \kappa 2\Big]\,\frac{g_0(x)}{\kappa}
-\frac{\rho(\kappa-4-\rho)}{2\kappa}.
\label{v}\EDE Since $R\in(0,1)$ on ${\cal D}'$, so $V_0(R)$ is well defined on ${\cal D}'$.
From (\ref{1+R})-(\ref{v}), we have that
 \begin{eqnarray}
\frac{\pa_j V_0(R)}{V_0(R)}&=&\frac{g_0(R)}\kappa\,\Big(\frac{A_{j,1}}{E_{j,j}}-\frac{A_{j,1}}{E_{j,k}}\Big)\,{\pa\xi_j(t_j)}
+\frac{g_0(R)}\kappa\,\Big(\frac{\kappa}2-3\Big)\Big[\Big(\frac{A_{j,2}}{E_{j,j}}-\frac{A_{j,2}}{E_{j,k}}\Big)\nonumber\\
&-&\Big(\frac{2A_{j,1}^2}{E_{j,0} E_{j,j}}-\frac{2A_{j,1}^2}{E_{j,0} E_{j,k}}\Big)\Big]\,{\pa t_j}
-\frac{\rho(\kappa-4-\rho)}{2\kappa}\,
\Big(\frac{A_{j,1}}{E_{j,k}}-\frac{A_{j,1}}{E_{j,j}}\Big)^2\,{\pa t_j}.\label{SDE-Y*}\end{eqnarray}

Define $N$ and $F$ on $\cal D$ such that $N=\frac{A_{1,1}A_{2,1}}{(A_{1,0}-A_{2,0})^2}$ and
$F(t_1,t_2)=\exp (\int_0^{t_2}\int_0^{t_1}2N(s_1,s_2)ds_1ds_2 )$. Let $\alpha
=\frac{6-\kappa}{2\kappa}$ and $\lambda=\frac{(8-3\kappa)(6-\kappa)}{2\kappa}$.
The following equations are (4.13) in \cite{reversibility} and (4.25) in \cite{duality}.
\begin{eqnarray} \frac{\pa_j N^\alpha}{N^\alpha}&=&\frac1\kappa\Big(3-\frac\kappa2\Big)\Big(\frac{A_{j,2}}{A_{j,1}}-\frac{2A_{j,1}}{E_{j,0}}\Big)\pa\xi_j(t_j)
+\lambda\Big(\frac 14\cdot\frac{A_{1,2}^2}{A_{1,1}^2}-\frac16\cdot\frac{A_{1,3}}{A_{1,1}}\Big)\,\pa t_j;\label{SDE-N*}\\
\frac{\pa_j F^{-\lambda}}{F^{-\lambda}}&=&-\lambda\Big(\frac
14\cdot \frac{A_{j,2}^2}{A_{j,1}^2}-\frac 16\cdot\frac{A_{j,3}}
{A_{j,1}}\Big)\,\pa t_j.\label{E*}\end{eqnarray}

Let $\tau=\frac{(\rho+2)(\kappa-6-\rho)}{2\kappa}$ and $\delta=-\frac{\rho(\kappa-4-\rho)}{4\kappa}$. Define $M$ on ${\cal D}'$ such that
\BGE M=|x_1-x_2|^{\tau}r_1(t_1)r_2(t_2)
D ^{\delta}V_0(R) N ^{\alpha}F ^{-\lambda}.\label{M*}\EDE
From (\ref{kappa-rho-deg*-1}), (\ref{ODE-r*}), (\ref{ODE-D*}) and (\ref{SDE-Y*})-(\ref{E*}),    we get
$$\frac{\pa_j   M}{  M}=\Big[\Big(3-\frac\kappa2\Big)\Big(\frac{A_{j,2}}{A_{j,1}}-\frac{2A_{j,1}}{E_{j,0}}\Big)+
g_0(R)\Big(\frac{A_{j,1}}{E_{j,j}}-\frac{A_{j,1}}{E_{j,k}}\Big)$$
\BGE-\frac{\rho}{\xi_j(t_j)-p_j(t_j)}-\frac{{\kappa-6-\rho} }{\xi_j(t_j)-q_j(t_j)}\Big]\,
\frac{dB_j(t_j)}{\sqrt\kappa}.\label{SDE-M*}\EDE

Define $\til r_j$ on $[0,T_j)$ such that \BGE\til r_j(t_j)=|\xi_j(t_j)-q_j(t_j)|^{-\frac{\kappa-6-\rho}{\kappa}}
|q_j(t_j)-p_j(t_j)|^{-\frac{\rho(\kappa-6-\rho)}{2\kappa}}\pa_z \vphi_j(t_j,x_{3-j})^{\frac{(\rho+2)(\kappa-6-\rho)}{4\kappa}}.\label{til-r}\EDE
Define $\til M$ on $\cal D$ such that
\BGE \til M=|x_1-x_2|^{\tau}{\til r_1(t_1) \til r_2(t_2)  D ^{\delta}|E_{1,2}E_{2,1}|^{-\frac\rho\kappa}U_0
(R) N^\alpha F ^{-\lambda}}.\label{Mc*}\EDE
Then $\til M$ is continuous on $\cal D$. Define $L_j$ on $\cal D$ such that if $t_j=0$ then $L_j=\pa_z \vphi_k(t_k,x_j)$;  if  $t_j>0$ then
\BGE L_j(t_1,t_2)=\frac{|E_{j,j}(t_1,t_2)|}{|\xi_j(t_j)-p_j(t_j)|}=
\frac{\vphi_{k,t_j}(t_k,\xi_j(t_j))-\vphi_{k,t_j}(t_k,p_j(t_j))}{\xi_j(t_j)-p_j(t_j)}.\label{Hj*}\EDE
Here the second ``$=$'' holds because $E_{j,j}$ has the same sign as ${\xi_j(t_j)-p_j(t_j)}$.
Since $\lim_{t_k\to 0^+}\xi_k(t_k)=\lim_{t_k\to 0^+} p_k(t_k)=x_k$ and $\lim_{t_j\to 0^+}\vphi_{k,t_j}(t_k,\cdot)
=\vphi_{k,0}(t_k,\cdot)=\vphi_k(t_k,\cdot)$, so
$L_j$ is continuous on $\cal D$. From (\ref{r*}), (\ref{M*}), (\ref{til-r})-(\ref{Hj*}), and that $V_0(x)=x^{\frac\rho\kappa}U_0(x)$, we find that
$M=\til M L_1^{\frac\rho\kappa}L_2^{\frac\rho\kappa}$ on ${\cal D}'$. Thus $M$ has continuous extension
to $\cal D$. Now we check the value of $M$ when $t_j=0$.

We have $\xi_j(0)=p_j(0)=x_j$, $q_j(0)=x_k$, and $K_j(0)=\emptyset$.
So $K_j(0)\cup K_k(t_k)=K_k(t_k)$. From (\ref{K-j-t*}) we have $K_{k,0}(t_k)=K_k(t_k)$ and $K_{j,t_k}(0)=\emptyset$, which implies that $\vphi_{k,0}(t_k,\cdot)=\vphi_k(t_k,\cdot)$ and $\vphi_{j,t_k}(0,\cdot)=\id$. Thus, if $t_j=0$, then $\til r_j(t_j)=|x_j-x_k|^{-\tau}$; and
$A_{j,0}=\vphi_k(t_k,x_j)=q_k(t_k)=B_{j,0}$,  $A_{j,1}=\pa_z\vphi_k(t_k,x_j)=\til B_{j,1}$, $A_{j,2}=\pa_z^2\vphi_k(t_k,x_j)$,
 $A_{k,0}=\xi_k(t_k)$, $B_{k,0}=p_k(t_k)$, and $A_{k,1}=1=\til B_{k,1}$, which imply that
 $E_{j,j}=0$, $E_{j,k}=q_k(t_k)-p_k(t_k)$, $E_{k,0}=E_{k,j}=\xi_k(t_k)-q_k(t_k)=-E_{j,0}$, $E_{k,k}=\xi_k(t_k)-p_k(t_k)$,
$|C_{j,k}|=|p_k(t_k)-q_k(t_k)|$, $D=\frac{\pa_z\vphi_k(t_k,x_j)}{(p_k(t_k)-q_k(t_k))^2}$, $R=0$, $U_0(R)=1$,
$N=\frac{\pa_z\vphi_k(t_k,x_j)}{(\xi_k(t_k)-q_k(t_k))^2}$, and $F=1$.
From (\ref{til-r}), (\ref{Mc*}) and the above argument, we find that
$\til M=\pa_z \vphi_k(t_k,x_j)^{-\frac\rho\kappa}$  when $t_j=0$. From the definition, $L_j=\pa_z \vphi_k(t_k,x_j)$ when $t_j=0$.
Since $\vphi_{j,t_k}(0,\cdot)=\id$, so $L_k=1$ when $t_j=0$. Thus, after continuous extension,
 $M=1$ when $t_1$ or $t_2$ equals $0$.

Let $Q_j$ be the formula inside the square bracket in (\ref{SDE-M*}), that is,
\BGE Q_j=\Big(3-\frac\kappa2\Big)\Big(\frac{A_{j,2}}{A_{j,1}}-\frac{2A_{j,1}}{E_{j,0}}\Big)+
g_0(R)\Big(\frac{A_{j,1}}{E_{j,j}}-\frac{A_{j,1}}{E_{j,k}}\Big)-\frac{\rho}{\xi_j(t_j)-p_j(t_j)}-\frac{{\kappa-6-\rho} }{\xi_j(t_j)-q_j(t_j)}.
\label{def-Q}\EDE
Then $Q_j$ is defined on ${\cal D}'$. Using the  observation in the previous paragraph and the fact that $g_0(0)=\rho$ and
$g_0$ is differentiable at $0$, we may check that $Q_j$ has continuous extension to $\cal D$. Thus, after   continuous extensions, the formula
$\frac{\pa_j M}{M}=Q_j\frac{dB_j(t_j)}{\sqrt\kappa}$ holds in $\cal D$.
For each $t_k\in[0,T_k)$, let $T_j(t_k)$ be the maximal number such that $K_j(t)\cap K_k(t_k)=\emptyset$ for $0\le t<T_j(t_k)$.
From (\ref{SDE-M*}) we conclude that for any fixed stopping time $t_k\in[0,T_k)$, $M$ is a continuous local martingale in $t_j$, where
$t_j$ ranges in $[0,T_j)$.

 Let $\HP$ denote the set of $(H_1,H_2)$ such that for $j=1,2$, $H_j$ is a hull in $\HH$ w.r.t.\ $\infty$ that contains some neighborhood of $x_j$ in
$\HH$, and $\lin{H_1}\cap\lin{H_2}=\emptyset$. For $(H_1,H_2)\in\HP$, let $T_j(H_j)$ be the first $t$ such that $\beta_j(t_j)
\in\lin{\HH\sem H_j}$. Then $T_j(H_j)$ is an $(\F^j_t)$-stopping time.

\begin{Theorem} For any $(H_1,H_2)\in\HP$, there are $C_2>C_1>0$ depending on $H_1$ and $H_2$
 such that $C_1\le M(t_1,t_2)\le C_2$ for $(t_1,t_2)\in[0,T_1(H_1)]\times [0,T_2(H_2)]$. \label{bound}
\end{Theorem}
{\bf Proof.} 
Since
$M=\til M L_1^{\frac\rho\kappa}L_2^{\frac\rho\kappa}$,  so we suffice to show that the theorem holds for $\til M$ and $L_j$, $j=1,2$.
To check the boundedness of $\til M$, we suffice to show that the theorem holds for every factor on the right-hand side of (\ref{Mc*}).
From Lemma \ref{f>*}, we find that the theorem holds for $U_0(R)$. The boundedness of other factors in (\ref{Mc*}) can be proved using
the method in Section 5 of \cite{reversibility}. For the boundedness of $L_j$, we suffice to note that from Lemma 5.2 in
\cite{reversibility}, the value of $L_j$ lies between $A_{j,1}$ and $\til B_{j,1}$, which are both uniformly bounded from $\infty$
and $0$. $\Box$

\vskip 3mm

Fix $(H_1,H_2)\in\HP$. From the local martingale property of $M$ and the above theorem, we see that $\EE[M(T_1(H_1),T_2(H_2))]=1$.
Let $\mu$ denote the joint distribution of $(\xi_1(t),0\le t<T_1)$ and $(\xi_2(t),0\le t<T_2)$. Define $\nu$ such that
$d\nu/d\mu=M(T_1(H_1),T_2(H_2))$. Then $\nu$ is also a probability measure. Suppose temporarily that the joint distribution
of $\xi_1$ and $\xi_2$ is $\nu$ instead of $\mu$. For $(t_1,t_2)\in\cal D$, define
\BGE B_{1,t_2}(t_1)=B_1(t_1)-\frac 1{\sqrt\kappa}\int_0^{t_1} Q_1(s,t_2)ds,\quad B_{2,t_1}(t_2)=B_2(t_2)-
\frac 1{\sqrt\kappa}\int_0^{t_2} Q_2(t_1,s)ds\label{Bjtk}.\EDE
Fix an $(\F^k_t)$-stopping time $\bar t_k$ with $\bar t_k\le T_k(H_k)$.
Since $B_j(t)$ is an $(\F^j_{t}\times \F^k_{\bar t_k})_{t\ge 0}$-Brownian motion under $\mu$, so
from (\ref{SDE-M*}), (\ref{def-Q}) and the Girsanov's Theorem, $B_{j,\bar t_k}(t)$, $0\le t\le T_j(H_j)$,
 is a partial $(\F^j_{t}\times \F^k_{\bar t_k})_{t\ge 0}$-Brownian motion under $\nu$.

\vskip 3mm

The following theorem is Theorem 6.1 in \cite{reversibility} and Theorem 4.5 in \cite{duality}. It
can be proved using the above theorem and the argument in \cite{reversibility} or \cite{duality}.

\begin{Theorem} For any $(H_1^m,H_2^m)\in\HP$, $1\le m\le n$, there
is a continuous function $M_*(t_1,t_2)$ defined on $[0,\infty]^2$
that satisfies the following properties: (i) $M_*=M$ on
$[0,T_1(H_1^m)]\times[0,T_2(H_2^m)]$ for $m=1,\dots,n$; (ii)
$M_*(t,0)=M_*(0,t)=1$ for any $t\ge 0$; (iii) $M_*(t_1,t_2)\in
[C_1,C_2]$ for any $t_1,t_2\ge 0$, where $C_2>C_1>0$ are constants
depending only on $H_j^m$, $j=1,2$, $1\le m\le n$; (iv) for any
$(\F^2_t)$-stopping time $\bar t_2$, $(M_*(t_1,\bar t_2),t_1\ge 0)$
is a bounded continuous $(\F^1_{t_1}\times\F^2_{\bar t_2})_{t_1\ge
0}$-martingale; and (v) for any $(\F^1_t)$-stopping time $\bar t_1$,
$(M_*(\bar t_1, t_2),t_2\ge 0)$ is a bounded continuous $(\F^1_{\bar
t_1}\times\F^2_{t_2})_{t_2\ge 0}$-martingale.
 \label{martg}
\end{Theorem}

\section{Coupling Measures}

 Let ${\cal C}:=\cup_{T\in(0,\infty]}
C([0,T))$. The map $T:{\cal C}\to (0,\infty]$ is such that
$[0,T(\xi))$ is the definition domain of $\xi$. For
$t\in[0,\infty)$, let $\F_t$ be the $\sigma$-algebra on $\cal C$
generated by $\{T>s,\xi(s)\in A\}$, where  $s\in[0,t]$ and $A$ is a Borel set on $\R$.
Then $(\F_t)$ is a
filtration on $\cal C$, and $T$ is an $(\F_t)$-stopping time. Let $\F_\infty=\vee_t \F_t$.

For $\xi\in\cal C$, let $K_\xi(t)$, $0\le t<T(\xi)$, denote the
chordal Loewner hulls driven by $\xi$.
Let $H$ be a hull in $\HH$ w.r.t.\ $\infty$. Let
$T_H(\xi)\in[0,T(\xi)]$ be the maximal number such that
$K_\xi(t)\cap\lin{\HH\sem H}=\emptyset$ for $0\le t<T_H$. Then $T_H$
is an $(\F_t)$-stopping time. Let ${\cal C}_H=\{T_H>0\}$. Then
$\xi\in{\cal C}_H$ iff $H$ contains some neighborhood of $\xi(0)$ in
$\HH$. Define $P_H:{\cal C}_H\to {\cal C}$ such that $P_H(\xi)$ is
the restriction of $\xi$ to $[0,T_H(\xi))$. Then $P_H({\cal
C}_H)=\{T_H=T\}$, and $P_H\circ P_H=P_H$.
If $A$ is a Borel set on $\R$ and $s\in[0,\infty)$, then
$$P_H^{-1}(\{\xi\in{\cal C}:T(\xi)>s,\xi(s)\in A\})=\{\xi\in{\cal
C}_H:T_H(\xi)>s,\xi(s)\in A\}\in\F_{T_H^-}.$$ Thus, $P_H$ is
$(\F_{T_H^-},\F_{\infty})$-measurable on ${\cal C}_H$. On the other
hand, the restriction of $\F_{T_H^-}$ to ${\cal C}_H$ is the
$\sigma$-algebra generated by $\{\xi\in{\cal
C}_H:T_H(\xi)>s,\xi(s)\in A\}$, where $s\in[0,\infty)$ and $A$ is a
Borel set on $\R$. Thus, $P_H^{-1}(\F_\infty)$ agrees with the
restriction of $\F_{T_H^-}$ to ${\cal C}_H$.

Let $\ha\C=\C\cup\{\infty\}$ be the Riemann sphere with spherical
metric. Let $\Gamma_{\ha\C}$ denote the space of nonempty compact
subsets of $\ha\C$ endowed with Hausdorff metric. Then
$\Gamma_{\ha\C}$ is a compact metric space. Define $G:{\cal C}\to
\Gamma_{\ha\C}$ such that $G(\xi)$ is the spherical closure of
$\{t+i\xi(t):0\le t<T(\xi)\}$. Then $G$ is a one-to-one map. Let
$I_G=G({\cal C})$. Let $\F^H_{I_G}$ denote the $\sigma$-algebra on
$I_G$ generated by Hausdorff metric. Let $${\cal
R}=\{\{z\in\C:a<\Ree z<b,c<\Imm z< d\}:a,b,c,d\in\R\}.$$ Then
$\F^H_{I_G}$ agrees with the $\sigma$-algebra on $I_G$ generated by
$\{\{F\in I_G:F\cap R\ne\emptyset\}:R\in{\cal R}\}$. Using this
result, one may check that $G$ and $G^{-1}$ (defined on $I_G$) are
both measurable with respect to $\F_\infty$ and $\F^H_{I_G}$.

Now we adopt the notation in the previous section.
Let $\mu_j$ denote the distribution of $(\xi_j(t),0\le t<T_j)$, which is a probability measure
on $\cal C$. Let $\mu=\mu_1\times\mu_2$ be a probability measure on
${\cal C}^2$. Since $\xi_1$ and $\xi_2$ are independent, so $\mu$ is the joint distribution of $\xi_1$ and
$\xi_2$.

Let $\HP_*$ be the set of $(H_1,H_2)\in\HP$ such that for $j=1,2$,
$H_j$ is a polygon whose vertices have rational coordinates. Then
$\HP_*$ is countable. Let $(H_1^m,H_2^m)$, $m\in\N$, be an
enumeration of $\HP_*$. For each $n\in\N$, let $M_*^n(t_1,t_2)$ be
the $M_*(t_1,t_2)$ given by Theorem \ref{martg} for $(H_1^m,H_2^m)$,
$1\le m\le n$, in the above enumeration. For each $n\in\N$ define
$\nu^n=(\nu^n_1,\nu^n_2)$ such that
${d\nu^n}/{d\mu}=M_*^n(\infty,\infty)$. From Theorem \ref{martg},
$M_*^n(\infty,\infty)>0$ and $\int
M_*^n(\infty,\infty)d\mu=\EE_\mu[M_*^n(\infty,\infty)]=1$, so
$\nu^n$ is a probability measure on ${\cal C}^2$. Since
$d\nu^n_1/d\mu_1=\EE_\mu[M_*^n(\infty,\infty)|\F^1_\infty]=M_*^n(\infty,0)=1$,
so $\nu^n_1=\mu_1$. Similarly, $\nu^n_2=\mu_2$. So each $\nu^n$ is a
coupling of $\mu_1$ and $\mu_2$.

Let $\bar\nu^n=(G\times G)_*(\nu^n)$ be a probability measure on
$\Gamma_{\ha \C}^2$. Since $\Gamma_{\ha\C}^2$ is compact, so
$(\bar\nu^n)$ has a subsequence $(\bar\nu^{n_k})$ that converges
weakly to some probability measure $\bar\nu=(\bar\nu_1,\bar\nu_2)$
on $\Gamma_{\ha\C}\times \Gamma_{\ha\C}$. Then for $j=1,2$,
$\bar\nu^{n_k}_j\to\bar\nu_j$ weakly. For $n\in\N$ and $j=1,2$,
since $\nu^n_j=\mu_j$, so $\bar \nu^n_j=G_*(\mu_j)$. Thus
$\bar\nu_j=G_*(\mu_j)$, $j=1,2$. So $\bar\nu$ is supported by
$I_G^2$. Let $\nu=(\nu_1,\nu_2)=(G^{-1}\times G^{-1})_*(\bar\nu)$ be
a probability measure on ${\cal C}^2$. Here we use the fact that
$G^{-1}$ is $(\F^H_{I_G},\F^j_\infty)$-measurable. For $j=1,2$, we
have $\nu_j=(G^{-1})_*(\bar\nu_j)=\mu_j$. So $\nu$ is also a
coupling measure of $\mu_1$ and $\mu_2$.

\vskip 3mm

The following lemma is Lemma 4.1 in \cite{duality}. The proof is similar.

\begin{Lemma} For any $n\in\N$, the restriction of $\nu$ to
$\F^1_{T_{H_1^n}}\times \F^2_{T_{H_2^n}}$ is absolutely continuous
w.r.t.\ $\mu$, and the Radon-Nikodym derivative is
$M(T_{H_1^n}(\xi_1),T_{H_2^n}(\xi_2))$. \label{Radon}
\end{Lemma}

Now suppose that the joint distribution of $\xi_1(t)$, $0\le t<T_1$, and
$\xi_2(t)$, $0\le t<T_2$, is the $\nu$ in the above lemma instead of $\mu=\mu_1\times\mu_2$.
Since the distribution of $\xi_j$ is $\nu_j=\mu_j$, so $\beta_j(t)$, $0\le t<T_j$,
is still a chordal SLE$(\kappa;\rho,\kappa-6-\rho)$ trace started from $(x_j;x_j^{\sigma_j},x_k)$.
Thus, a.s.\ $\lim_{t\to T_j^-} \beta_j(t)=x_k$. For $(t_1,t_2)\in\cal D$,
let $B_{j,t_k}(t_j)$ be defined by (\ref{Bjtk}). Fix an $(\F^k_t)$-stopping
time $\bar t_k\in[0,T_k)$. Choose any $n\in\N$. Let $\bar t^n_k=\bar t_k\wedge T_k(H^n_k)$.
Then $\bar t_k^n$ is also an $(\F^k_t)$-stopping time, and satisfies $\bar t_k^n\le T_k(H^n_k)$.
From the above lemma and the discussion after Theorem \ref{bound}, we see that
$B_{j,\bar t_k^n}(t)$, $0\le t\le T_j(H^n_j)$, is a partial $(\F^j_t\times \F^k_{\bar t_k^n})_{t\ge 0}$-Brownian motion.

\begin{Lemma} $B_{j,\bar t_k}(t)$, $0\le t<T_j(\bar t_k)$, is a partial  $(\F^j_t\times \F^k_{\bar t_k})_{t\ge 0}$-Brownian motion.
\label{BM}
\end{Lemma}
{\bf Proof.} Write $T_j^n$ for $T_j(H^n_j)$, $j=1,2$, $n\in\N$.
 From the above argument, we know that for any $n\in\N$, $B_{j,\bar t_k^n}(t\wedge T_j^n)$,
$0\le t<\infty$, is a continuous $(\F^j_t\times \F^k_{\bar t_k^n})_{t\ge 0}$-martingale. Define
$S^n_j=T^n_j$ on $\{\bar t_k\le T^n_k\}$, and $S^n_j=0$ on $\{T^n_k<\bar t_k\}$. Then for
any $t\ge 0$, $\{S^n_j\le t\}=\{T^n_k<\bar t_k\}\cup\{T^n_j\le t\}\in\F^j_t\times\F^k_{\bar t_k}$.
So $S^n_j$ is an $(\F^j_t\times \F^k_{\bar t_k})_{t\ge 0}$-stopping time. Now we claim that
 $B_{j,\bar t_k}(t\wedge S_j^n)$, $0\le t<\infty$, is a continuous $(\F^j_t\times \F^k_{\bar t_k})_{t\ge 0}$-martingale.
Fix $s_2\ge s_1\ge 0$ and ${\cal E}\in \F^j_{s_1}\times \F^k_{\bar t_k}$. Let ${\cal E}_1={\cal E}\cap\{T^n_k<\bar t_k\}$
and ${\cal E}_2={\cal E}\cap\{\bar t_k\le T^n_k\} $. Since $S_j^n=0$ on ${\cal E}_1$, so $B_{j,\bar t_k}(s_2\wedge S_j^n)=
0=B_{j,\bar t_k}(s_1\wedge S_j^n)$ on ${\cal E}_1$, which implies that
\BGE\int_{{\cal E}_1} B_{j,\bar t_k}(s_2\wedge S_j^n)\,d\nu=0=\int_{{\cal E}_1} B_{j,\bar t_k}(s_1\wedge S_j^n)\,d\nu.\label{mart-1}\EDE
Since $\bar t_k=\bar t_k^n$ on $\{\bar t_k\le T^n_k\}$, so $\F^k_{\bar t_k}$ agrees with $\F^k_{\bar t_k^n}$ on $\{\bar t_k\le T^n_k\}$. Thus,
${\cal E}_2\in \F^j_{s_1}\times \F^k_{\bar t_k^n}$. Since $\bar t_k=\bar t_k^n$ and $S_j^n=T_j^n$ on ${\cal E}_2$, so from the
martingale property of $B_{j,\bar t_k^n}(t\wedge T_j^n)$, we have
\BGE\int_{{\cal E}_2} B_{j,\bar t_k}(s_2\wedge S_j^n)\,d\nu=\int_{{\cal E}_2} B_{j,\bar t_k}(s_1\wedge S_j^n)\,d\nu.\label{mart-2}\EDE
Since $\cal E$ is the disjoint union of ${\cal E}_1$ and ${\cal E}_2$, so
from (\ref{mart-1}) and (\ref{mart-2}), $\EE_\nu[B_{j,\bar t_k}(s_2\wedge S_j^n)|\F^j_{s_1}\times \F^k_{\bar t_k}]
=B_{j,\bar t_k}(s_1\wedge S_j^n)$. So the claim is justified.

Since the above claim holds for any $n\in\N$, so $B_{j,\bar t_k}(t)$, $0\le t<\vee_{n=1}^\infty S_j^n$, is
 a continuous $(\F^j_t\times \F^k_{\bar t_k})_{t\ge 0}$-local martingale. We now claim that $\vee_{n=1}^\infty S_j^n=T_j(\bar t_k)$.
Fix any $n\in\N$. If $T^n_k<\bar t_k$ then $S_j^n=0< T_j(\bar t_k)$. If $\bar t_k\le T^n_k$ then $S_j^n=T_j^n$.
From $\bar t_k\le T^n_k$ we have $K_k({\bar t_k})\subset H_k^n$. From $S_j^n=T_j^n$ we have
 $K_j({S_j^n})\subset H_j^n$. Since $\lin{H_j^n}\cap \lin{H_k^n}=\emptyset$, so
  $\lin{K_j({S_j^n})}\cap \lin{K_k({\bar t_k})}=\emptyset$, and so again we have $S_j^n< T_j(\bar t_k)$. Since the above holds
  for any $n\in\N$, so $\vee_{n=1}^\infty S_j^n\le T_j(\bar t_k)$. Now suppose $t_0<T_j(\bar t_k)$. Then $\lin{K_j({t_0})}\cap \lin{K_k(\bar t_k)}
=\emptyset$. We may always find $(H^{n_0}_1,H^{n_0}_2)\in\HP_*$ such that $K_j(t_0)\subset H^{n_0}_j$ and $K_k(\bar t_k)\subset H^{n_0}_k$.
Then we have $\bar t_k\le T_k^{n_0}$. So $\vee_{n=1}^\infty S_j^n\ge S_j^{n_0}=T_j^{n_0}\ge t_0$. Since this holds for any $t_0<T_j(\bar t_k)$,
so $\vee_{n=1}^\infty S_j^n=T_j(\bar t_k)$. Thus, $B_{j,\bar t_k}(t)$, $0\le t<T_j(\bar t_k)$, is
 a continuous $(\F^j_t\times \F^k_{\bar t_k})_{t\ge 0}$-local martingale. Using a similar argument, we conclude that
 $B_{j,\bar t_k}(t)^2-t$, $0\le t<T_j(\bar t_k)$, is also
 a continuous $(\F^j_t\times \F^k_{\bar t_k})_{t\ge 0}$-local martingale.
Using the characterization of Brownian motion in \cite{RY}, we complete the proof. $\Box$

\begin{Theorem} Let $a>0$. Let $\bar t_2\in(0,T_2)$ be an $(\F^2_t)$-stopping time. Let
$C_1=a\cdot\frac{\xi_2(\bar t_2)-p_2(\bar t_2)}{p_2(\bar t_2)-q_2(\bar t_2)}>0$,
$w(z)=C_1\cdot\frac{z-q_2(\bar t_2)}{\xi_2(\bar t_2)-z}$, and $W=w\circ\vphi_2(\bar t_2,\cdot)$.
Then after a time-change, $W(\beta_1(t))$, $0\le t<T_1(\bar t_2)$, has the distribution of a degenerate  intermediate
SLE$(\kappa;\rho)$ trace with force points $0_+$ and $a$. Moreover, a.s.\ $T_1(\bar t_2)<T_1$ and
$\beta_1(T_1(\bar t_2))=\beta_2(\bar t_2)$.
\label{theorem-intermediate}
\end{Theorem}
{\bf Proof.} Let $C_2=C_1\cdot(\xi_2(\bar t_2)-q_2(\bar t_2))>0$.  For $0\le t<T_1(\bar t_2)$,
define
\BGE \til\vphi(t,z)=\frac{C_2A_{2,1}(t,\bar t_2)}{A_{2,0}(t,\bar t_2)-\vphi_{1,\bar t_2}(t,w^{-1}(z))}-C_1+
\int_0^t \frac{2C_2A_{2,1}(s,\bar t_2)A_{1,1}(s,\bar t_2)^2}{E_{1,0}(s,\bar t_2)^3}\,ds;\label{til-phi}\EDE
\BGE \til\xi(t)=\frac{C_2A_{2,1}(t,\bar t_2)}{E_{2,0}(s,\bar t_2)}-C_1+
\int_0^t \frac{2C_2A_{2,1}(s,\bar t_2)A_{1,1}(s,\bar t_2)^2}{E_{1,0}(s,\bar t_2)^3}\,ds;\label{til-xi}\EDE
\BGE\til p(t)=\frac{C_2 A_{2,1}(t,\bar t_2)}{E_{2,1}(t,\bar t_2)}-C_1+
\int_0^t \frac{2C_2A_{2,1}(s,\bar t_2)A_{1,1}(s,\bar t_2)^2}{E_{1,0}(s,\bar t_2)^3}\,ds;\label{til-p}\EDE
\BGE\til q(t)=\frac{C_2A_{2,1}(t,\bar t_2)}{E_{2,2}(t,\bar t_2)}-C_1+
\int_0^t \frac{2C_2A_{2,1}(s,\bar t_2)A_{1,1}(s,\bar t_2)^2}{E_{1,0}(s,\bar t_2)^3}\,ds.\label{til-q}\EDE

Since $A_{2,0}(0,\bar t_2)=\xi_2(\bar t_2)$,  $A_{2,1}(0,\bar t_2)=1$, and $\vphi_{1,\bar t_2}(0,\cdot)=\id$, so $\til\vphi(0,z)=z$.
Using (\ref{pa-A_2}) and (\ref{chordal*-}) with $j=1$ and $k=2$, it is straightforward to check that
\BGE \pa_t \til\vphi(t,z)= \frac{2C_2^2N(t,\bar t_2)^2}{\til\vphi(t,z)-\til\xi(t)}.\label{til-chordal}\EDE
Let $v(t)=\int_0^t C_2^2N(s,\bar t_2)^2ds$. Then $v(0)=0$ and $v$ is continuous and
strictly increasing. So $v$ maps $[0,T_1(\bar t_2))$ onto $[0,T)$ for some $T\in(0,\infty]$. Let $\vphi(t,\cdot)
=\til\vphi(v^{-1}(t),\cdot)$ and $\xi(t)=\til\xi(v^{-1}(t))$ %
 for $0\le t<T$. From (\ref{til-chordal}), we have
$\pa_t \vphi(t,z)=\frac{2}{\vphi(t,z)-\xi(t)}$. Thus $\vphi(t,\cdot)$, $0\le t<T$, are the chordal Loewner
maps driven by $\xi$.

Note that $w$ maps $\HH$ conformally onto $\HH$, and $w(\xi_2(\bar t_2))=\infty$. Since $\vphi_2(\bar t_2,\cdot)$
maps $\HH\sem \beta_2((0,\bar t_2])$ conformally onto $\HH$, and $\vphi_2(\bar t_2,\beta_2(\bar t_2))=\xi_2(\bar t_2)$, so
$W$ maps $\HH\sem \beta_2((0,\bar t_2])$ conformally on $\HH$, and $W(\beta_2(\bar t_2))=\infty$.
For any $t\in[0,T_1(\bar t_2))$, $w^{-1}$ maps $\HH\sem W(\beta_1((0,t]))$ conformally onto
$\HH\sem \vphi_2(\bar t_2,\beta_1((0,t]))=\HH\sem K_{1,\bar t_2}(t)$. Since $\vphi_{1,\bar t_2}(t,\cdot)$ maps
$\HH\sem K_{1,\bar t_2}(t)$ conformally onto $\HH$, so from (\ref{til-phi}),
$\til\vphi(t,\cdot)$ maps $\HH\sem W(\beta_1((0,t]))$ conformally onto $\HH$. For $0\le t<T$,
let $\beta(t)=W(\beta_1(v^{-1}(t)))$,
then $\vphi(t,\cdot)$ maps $\HH\sem \beta((0,t])$ conformally onto $\HH$. So $\beta(t)$, $0\le t<T$, is the chordal Loewner
 trace driven by $\xi$.

 Let $p(t)=\til p(v^{-1}(t))$ and $q(t)=\til q(v^{-1}(t))$, $0\le t<T$. Applying (\ref{pa-A_2}) and
(\ref{Bj0-1*}) with $j=1$ and $k=2$, and using $v'(t)=C_2^2N(t,\bar t_2)^2$, it is straightforward to check that
\BGE p'(t)=\frac{2}{p(t)-\xi(t)},\quad 0<t<T;\qquad q'(t)=\frac{2}{q(t)-\xi(t)}, \quad 0\le t\le T.\label{pq}\EDE
Moreover, since $A_{1,0}(t,\bar t_2)<B_{1,0}(t,\bar t_2)<B_{2,0}(t,\bar t_2)<A_{2,0}(t,\bar t_2)$ for $0< t<T_1(\bar t_2)$,
so from (\ref{til-xi})-(\ref{til-q}) and the definition of $E_{2,m}$, $m=0,1,2$, we have
\BGE \xi(t)<p(t)<q(t)<\infty,\quad 0< t<T.\qquad \label{xi<p<q}\EDE
Since $A_{1,0}(0,\bar t_2)=q_2(\bar t_2)=B_{1,0}(0,\bar t_2)$, and $A_{2,0}(0,\bar t_2)=\xi_2(\bar t_2)$, so
$E_{2,0}(0,\bar t_2)=E_{2,1}(0,\bar t_2)=\xi_2(\bar t_2)-q_2(\bar t_2)$. Note that $A_{2,1}(0,\bar t_2)=1$, so
\BGE\xi(0)=p(0)=\frac{C_2}{\xi_2(\bar t_2)-q_2(\bar t_2)}
-C_1=0.\label{xi=p=0}\EDE
Since $B_{2,0}(0,\bar t_2)=p_2(\bar t_2)$, so $E_{2,2}(0,\bar t_2)=\xi_2(\bar t_2)-p_2(\bar t_2)$. Thus,
\BGE q(0)=\frac{C_2}{\xi_2(\bar t_2)-p_2(\bar t_2)}-C_1
=a>0.\label{q=1}\EDE

Note that $E_{2,0}=-E_{1,0}$.
Applying (\ref{pa-A_2}) and (\ref{Djm*}) with $j=1$, $k=2$ and $m=0$, we get
\BGE d\til \xi(t)=C_2{N(t,\bar t_2)}d\xi_1(t)+C_2\,\frac{A_{2,1}(t,\bar t_2)}{E_{1,0}(t,\bar t_2)}
\,\Big[\Big(\frac \kappa 2-3\Big)\frac{A_{1,2}(t,\bar t_2)}{E_{1,0}(t,\bar t_2)}+(6-\kappa)\frac{A_{1,1}(t,\bar t_2)^2}
{E_{1,0}(t,\bar t_2)^2}\Big]\,dt.\label{SDE-xi-1}\EDE
From (\ref{kappa-rho-deg*-1}),  (\ref{def-Q}) and (\ref{Bjtk}),
 we see that $\xi_1(t)$, $0\le t<T_1(\bar t_2)$, satisfies the  $(\F^1_t\times \F^2_{\bar t_2})_{t\ge 0}$-adapted SDE:
\BGE d\xi_1(t)=\sqrt\kappa d  B_{1,\bar t_2}(t)+\Big[\Big(3-\frac\kappa2\Big)\Big(\frac{A_{1,2}}
 {A_{1,1}}-\frac{2A_{1,1}}{E_{1,0}}\Big)+
g_0(R)\Big(\frac{A_{1,1}}{E_{1,1}}-\frac{A_{1,1}}{E_{1,2}}\Big)\Big]\Big|_{(t,\bar t_2)}\,dt.\label{SDE-xi-2}\EDE
From (\ref{SDE-xi-1}) and (\ref{SDE-xi-2}) we conclude that
\BGE d\til \xi(t)=C_2N(t,\bar t_2)\Big[\sqrt \kappa d  B_{1,\bar t_2}(t)
+g_0(R(t,\bar t_2))\Big(\frac{A_{1,1}(t,\bar t_2)}{E_{1,1}(t,\bar t_2)}-\frac{A_{1,1}(t,\bar t_2)}{E_{1,2}(t,\bar t_2)}\Big)\,dt\Big].
\label{til-xi-SDE}\EDE
Let
\BGE S(t)=\frac{g_0(R(t,\bar t_2))}{C_2N(t,\bar t_2)}
\,\Big(\frac{A_{1,1}(t,\bar t_2)}{E_{1,1}(t,\bar t_2)}-\frac{A_{1,1}(t,\bar t_2)}{E_{1,2}(t,\bar t_2)}\Big).\label{S}\EDE
Since $\til\xi(t)=\xi(v(t))$ and $v'(t)=C_2^2N(t,\bar t_2)^2$, so from (\ref{til-xi-SDE}) and
 Lemma \ref{BM}, there is a
Brownian motion $B(t)$ such that for $0< t<T$,
\BGE d\xi(t)=\sqrt\kappa dB(t)+S(v^{-1}(t))dt.\label{xi-SDE}\EDE
From (\ref{til-xi})-(\ref{til-q}), we have
$$\til p(t)-\til\xi(t)=C_2\,\frac{A_{2,1}(t,\bar t_2)E_{1,1}(t,\bar t_2)}{E_{1,0}(t,\bar t_2)E_{2,1}(t,\bar t_2)};$$
$$\til q(t)-\til\xi(t)=C_2\,\frac{A_{2,1}(t,\bar t_2)E_{1,2}(t,\bar t_2)}{E_{1,0}(t,\bar t_2)E_{2,2}(t,\bar t_2)}.$$
Thus,
$$\frac{\til p(t)-\til \xi(t)}{\til q(t)-\til\xi(t)}=\frac{E_{1,1}(t,\bar t_2)E_{2,2}(t,\bar t_2)}{E_{1,2}(t,\bar t_2)E_{2,1}(t,\bar t_2)}
=R(t,\bar t_2).$$
From (\ref{J}), (\ref{S}) and the above formulas, we get
$$J(\til p(t)-\til\xi(t),\til q(t)-\til\xi(t))=-\Big(
\frac{1}{\til p(t)-\til \xi(t)}-\frac{1}{\til q(t)-\til\xi(t)}\Big)\cdot g_0\Big(\frac{\til p(t)-\til \xi(t)}{\til q(t)-\til\xi(t)}\Big)
=S(t).$$
From (\ref{xi-SDE}) we find that, for $0<t<T$,
\BGE d\xi(t)=\sqrt\kappa dB(t)+J(p(t)-\xi(t), q(t)-\xi(t))dt.\label{xi-SDE-2}\EDE
So $\xi(t)$, $p(t)$ and $q(t)$, $0<t<T$, solve (\ref{pq}) and (\ref{xi-SDE-2}), and satisfy (\ref{xi<p<q}-\ref{q=1}).
Assume that this solution  can be extended beyond $T$. Since $\kappa\in(0,4)$, so
$\beta(T)=\lim_{t\to T^-}\beta(t)\in\HH$. Thus, $\lim_{t\to (T_1(\bar t_2))^-} W(\beta_1(t))\in\HH$. From the definition,
$W$ maps $\HH\sem \beta((0,\bar t_2])$ conformally onto $\HH$. So we have
$\lim_{t\to (T_1(\bar t_2))^-} \beta_1(t)\in\HH\sem \beta((0,\bar t_2])$. This implies that the distance between
$\beta_1((0,T_1(\bar t_2)])$ and $\beta_2((0,\bar t_2])$ is positive. This is impossible because of the definition
of $T_1(\bar t_2)$ and the fact that $\lim_{t\to T_1^-} \beta_1(t)=x_2=\beta_2(0)$. Thus $(0,T)$ is the maximal interval of the solution.
From (\ref{pq})-(\ref{q=1}) and (\ref{xi-SDE-2}), we see that $\beta(t)$, $0\le t<T$,
is a degenerate  intermediate SLE$(\kappa;\rho)$ trace with force points $0^+$ and $a$. Since $\beta$ is a time-change of
$W(\beta_1)$, so after a time-change, $W(\beta_1(t))$, $0\le t<T_1(\bar t_2)$, has the distribution of a degenerate  intermediate
SLE$(\kappa;\rho)$ trace with force points $0_+$ and $a$.

From Corollary \ref{T=infty-degenerate*} and the fact that $W^{-1}(\infty)=\beta_2(\bar t_2)$,
we see that a.s.\ $\beta_2(\bar t_2)$ is a subsequential limit of $\beta_1(t)$ as $t\to (T_1(\bar t_2))^-$. If $T_1(\bar t_2)=T_1$
then $\lim_{t\to (T_1(\bar t_2))^-}\beta_1(t)=\lim_{t\to T_1^-}\beta_1(t)=x_2\ne \beta_2(\bar t_2)$ because $\bar t_2>0$, which
a.s.\ does not happen. Thus, a.s.\ $T_1(\bar t_2)<T_1$. Since $\beta_1$ is continuous on $[0,T_1)$, so a.s.\
 $\beta_1(T_1(\bar t_2))=\lim_{t\to (T_1(\bar t_2))^-}\beta_1(t)$.
Since a.s.\ $\beta_2(\bar t_2)$ is a subsequential limit of $\beta_1(t)$ as $t\to (T_1(\bar t_2))^-$, so
$\beta_1(T_1(\bar t_2))=\beta_2(\bar t_2)$. $\Box$

\begin{Theorem} Almost surely $\beta_1((0,T_1))=\beta_2((0,T_2))$.
\label{overlap}
\end{Theorem}
{\bf Proof.} For $n\in\N$, let $S_n$ be the first time that $|\beta_2(t)-x_1|=|x_2-x_1|/(n+1)$. Then
for each $n\in\N$, $S_n$ is an $(\F^2_t)$-stopping time, $S_n\in(0,T_2)$, and $T_2=\vee_{n=1}^\infty S_n$.
For each $q\in\Q_{>0}$, let $S_{n,q}=S_n\wedge q$, which is also an $(\F^2_t)$-stopping time.
Then $\{S_{n,q}:n\in\N,q\in\Q_{>0}\}$ is a dense subset of $(0,T)$. Applying Theorem \ref{theorem-intermediate} with
$\bar t_2=S_{n,q}$, we see that a.s.\
$\beta_2(S_{n,q})\in\beta_1((0,T_1))$ for any $n\in\N$ and $q\in\Q_{>0}$. From the denseness of $\{S_{n,q}\}$ and
the continuity of $\beta_1$, we have a.s.\ $\beta_2((0,T_2))\subset\beta_1((0,T_1))$. Since both $\beta_1$ and
$\beta_2$ are simple curves,   $\beta_1(0)=x_1=\beta_2(T_2)$, and $\beta_2(0)=x_2=\beta_1(T_1)$, so a.s.\
  $\beta_1((0,T_1))=\beta_2((0,T_2))$. $\Box$

\begin{Corollary} Suppose $\beta(t)$, $0\le t<\infty$, is a degenerate  intermediate SLE$(\kappa;\rho)$ trace. Then
a.s.\ $\lim_{t\to\infty}\beta(t)=\infty$.
\label{lim=infty}
\end{Corollary}
{\bf Proof.} Suppose that the force points for $\beta$ is $0^+$ and $a_0>0$. Applying
Theorem \ref{theorem-intermediate} with $a=a_0$ and any $(\F^2_t)$-stopping time $\bar t_2
\in(0,T_2)$. Then
$W(\beta_1(t))$, $0\le t<T_1(\bar t_2)$, has the same distribution as $\beta(t)$, $0\le t<\infty$, up to a time-change,
and a.s.\ $\lim_{t\to (T_1(\bar t_2))^-} \beta_1(t)=\beta_1(T_1(\bar t_2))=\beta_2(\bar t_2)$. Since $W(\beta_2(\bar t_2))=\infty$,
so a.s.\ $\lim_{t\to (T_1(\bar t_2))^-} W(\beta_1(t))=\infty$. Thus, a.s.\ $\lim_{t\to \infty}\beta(t)=\infty$.  $\Box$

\vskip 4mm

\no {\bf Proof of Theorem \ref{reversal*}.} We may find $W_1$ that maps $\HH$ conformally or conjugate conformally onto $\HH$ such that $W_1(x_1)=0$,
$W_1(x_1^+)=0^\sigma$, and $W_1(x_2)=\infty$. Let $W_2=W_0^{-1}\circ W_1$. Then
 $W_2$ maps $\HH$ conjugate conformally or conformally onto $\HH$ such that $W_2(x_2)=0$, $W_2(x_2^-)=0^\sigma$, and
$W_2(x_1)=\infty$. Recall that for $j=1,2$, $\beta_j(t)$, $0< t<T_j$, is a chordal SLE$(\kappa;\rho,\kappa-6-\rho)$ trace started
from $(x_j;x_j^{\sigma_j},x_{3-j})$, where $\sigma_1=+$ and $\sigma_2=-$.
From Proposition \ref{coordinate}, after a time-change, $W_j^{-1}(\beta_0(t))$, $0<t<\infty$, has the same distribution as
$\beta_j(t)$, $0< t<T_j$, $j=1,2$. 
From Theorem \ref{theorem-intermediate}, after a time-change, the reversal of $\beta_2(t)$, $0<t<T_2$, agrees with $\beta_1(t)$, $0<t<T_1$.
Thus, $W_2^{-1}(\beta_0(1/t))$, $0< t<\infty$, has the same distribution as $W_1^{-1}(\beta_0(t))$, $0<t<\infty$, after a time-change.
Since $W_0=W_1\circ W_2^{-1}$, so the proof is finished. $\Box$

\vskip 4mm

\no {\bf Proof of Theorem \ref{reversal-2}.} Applying Theorem \ref{theorem-intermediate} with any $(\F^2_t)$-stopping time
$\bar t_2\in(0,T_2)$ and $a=1/b_0$, we get $w(z)=a\cdot\frac{\xi_2(\bar t_2)-p_2(\bar t_2)}{p_2(\bar t_2)-q_2(\bar t_2)}
\cdot\frac{z-q_2(\bar t_2)}{\xi_2(\bar t_2)-z}$ and
 $W=w\circ\vphi_2(\bar t_2,\cdot)$, such that after a time-change,
$W(\beta_1(t))$, $0\le t<T_1(\bar t_2)$, has the same distribution as a degenerate  intermediate SLE$(\kappa;\rho)$ trace
with force points $0^+$ and $a=1/b_0$.

Let $\til T=T_2-\bar t_2$. For $0\le t<\til T$, let $\til\xi(t)=\xi_2(\bar t_2+t)$, $\til p(t)=p_2(\bar t_2+t)$ and $\til q(t)=q_2(\bar t_2+t)$.
Let $\til B(t)=B_2(\bar t_2+t)-B_2(\bar t_2)$, $t\ge 0$. Then $\til B(t)$ is a Brownian motion that is independent of $\xi_2(\bar t_2)$,
$p_2(\bar t_2)$ and $q_2(\bar t_2)$. From
(\ref{kappa-rho-deg*-1})-(\ref{kappa-rho-deg*-3}), $\til\xi(t)$, $\til p(t)$ and $\til q(t)$, $0\le t<\til T$, satisfy the following SDE:
\begin{eqnarray}
d\til\xi(t)& = &\sqrt\kappa d\til B(t)+\frac{\rho}{\til\xi(t)-\til p(t)}\,dt+\frac{\kappa-6-\rho}
{\til\xi(t)-\til q(t)}\,dt,\nonumber\\
d\til p(t)&=&\frac{2}{\til p(t)-\til \xi(t)}\,dt,\qquad  d\til q(t)\;\;=\;\:\frac{2}{\til q(t)-\til\xi(t)}\,dt\nonumber,
\end{eqnarray}
with initial values $$\til\xi(0)=\xi_2(\bar t_2),\quad \til p(0)=p_2(\bar t_2),\quad\til q(0)=q_2(\bar t_2).$$
For $0\le t<\til T$, let $\til\vphi(t,\cdot)=\vphi_2(\bar t_2+t,\cdot)\circ\vphi_2(\bar t_2,\cdot)^{-1}$
and $\til\beta(t) =\vphi_2(\bar t_2,\beta_2(\bar t_2+t))$.
Then $\til\vphi(0,z)=z$, and $\til\vphi(t,z)$, $0\le t<\til T$, satisfy
$\pa_t \til\vphi(t,z)=\frac{2}{\til\vphi(t,z)-\til\xi(t)}$, and for each $0\le t<\til T$,
$\til\vphi(t,\cdot)$ maps $\HH\sem\til\beta((0,t])$ conformally onto $\HH$. Thus, $\til\beta(t)$, $0\le t<\til T$,
is the chordal Loewner trace driven by $\til\xi$.
The solution $\til \xi(t)$, $\til p(t)$ and $\til q(t)$, $0\le t<\til T$, could not
be extended beyond $\til T$ because $\lim_{t\to \til T^-}\til\beta(t)=\vphi_2(\bar t_2,\lim_{t\to T_2^-}\beta_2(t))
=\vphi_2(\bar t_2,x_1)\in\R$. Thus, $\til\beta(t)=\vphi_2(\bar t_2,\beta_2(\bar t_2+t))$, $0<t<T_2-\bar t_2$,
is a chordal SLE$(\kappa;\rho,\kappa-6-\rho)$ trace started from $(\xi_2(\bar t_2);p_2(\bar t_2),q_2(\bar t_2))$.
Let  $W_1=W_0^{-1}\circ w$. Then $W_0^{-1}\circ W=W_1\circ\vphi_2(\bar t_2,\cdot)$,  $W_1$ maps $\HH$ conformally onto $\HH$,
$W_1(\xi_2(\bar t_2))=0$, $W_1(q_2(\bar t_2))=\infty$ and $W_1(p_2(\bar t_2))=1/a=b_0$. From Proposition \ref{coordinate},
$W_0^{-1}\circ W ( \beta_2(\bar t_2+t))=W_1(\til\beta(t))$, $0<t<T_2-\bar t_2$, has the same distribution as $\beta_0(t)$, $0<t<\infty$,
after a time-change. From Theorem \ref{theorem-intermediate} and Theorem \ref{overlap}, after a time-change, the
reversal of $\beta_2(t)$, $\bar t_2<t<T_2$, has the same distribution as $\beta_1(t)$, $0< t<T_1(\bar t_2)$.
Thus, after a time-change,  $W_0(\beta_0(1/t))$, $0<t<\infty$, has the same distribution as the reversal of
$W(\beta_1(t))$, $0< t<T_1(\bar t_2)$, which has the same distribution as a degenerate  intermediate SLE$(\kappa;\rho)$ trace
with force points $0^+$ and $1/b_0$. $\Box$

\vskip 3mm

Now we will see some applications of Theorem \ref{reversal*}. The following proposition is   Theorem
5.4 in \cite{duality}, where $\pa_\HH^+ S$ is the right boundary of $S$ in $\HH$ (c.f.\ \cite{duality}).

\begin{Proposition} Let $\kappa>4$, $C\ge 1/2$,
 and $K(t)$, $0\le t<\infty$, be a
chordal SLE$(\kappa;C(\kappa-4))$ process started from $(0;0^+)$.
Let $K(\infty)=\cup_{t<\infty}K(t)$. Let $W_0(z)=1/\lin z$. Then
$W_0(\pa_\HH^+ K(\infty))$ has the same distribution as the image of a
chordal SLE$(\kappa';C'(\kappa'-4),\frac 12(\kappa'-4))$ trace
started from $(0;0^+,0^-)$, where $\kappa'=16/\kappa$ and $C'=1-C$.
\label{duality2}
\end{Proposition}

Applying the above proposition with $C=1$, and applying
 Theorem \ref{reversal*} with $\kappa=\kappa'$ and $\rho=\frac 12(\kappa'-4)$, we
conclude the following theorem, which is Conjecture 2 in \cite{Julien-Duality}.

\begin{Theorem} Let $\kappa>4$,  and $K(t)$, $0\le t<\infty$, be a
chordal SLE$(\kappa;\kappa-4)$ process started from $(0;0^+)$.
Let $K(\infty)=\cup_{t<\infty}K(t)$. Then
$\pa_\HH^+ K(\infty)$ has the same distribution as the image of a
chordal SLE$(\kappa'; \frac 12(\kappa'-4))$ trace
started from $(0;0^-)$, where $\kappa'=16/\kappa$.
\label{JD}
\end{Theorem}

The following proposition is a part of Theorem 5.2 in \cite{duality2}.

\begin{Proposition} Let $\kappa>4$ and $C_+,C_-\ge 1/2$. Let
$K(t)$, $0\le t<\infty$, be a
chordal SLE$(\kappa;C_+(\kappa-4),C_-(\kappa-4))$ process started
from $(0;0^+,0^-)$. Let $K(\infty)=\cup_{t\ge 0}K(t)$. Let
 $\kappa'=16/\kappa$ and $W_0(z)=1/\lin{z}$. Then $W_0(\pa_\HH^+ K(\infty))$ has the same distribution
as the image of a chordal SLE$(\kappa';(1-C_+)(\kappa'-4),(1/2-C_-)(\kappa'-4))$ trace started
from $(0;0^+,0^-)$.
\label{duality1}
\end{Proposition}

Applying Proposition \ref{duality1} with $C+=1$ or $C-=1/2$, and using Theorem \ref{reversal*}, we
conclude the following two theorems.

\begin{Theorem} Let $\kappa>4$, $C\ge 1/2$,
 and $K(t)$, $0\le t<\infty$, be a
chordal SLE$(\kappa;\kappa-4,C(\kappa-4))$ process started
from $(0;0^+,0^-)$. Let $K(\infty)=\cup_{t<\infty}K(t)$. Then $\pa_\HH^+ K(\infty)$ has the same distribution
as the image of a chordal SLE$(\kappa';(1/2-C)(\kappa'-4))$ trace started
from $(0;0^-)$, where $\kappa'=16/\kappa$.
\label{duality1*}
\end{Theorem}

\begin{Theorem} Let $\kappa>4$, $C\ge 1/2$,
 and $K(t)$, $0\le t<\infty$, be a
chordal SLE$(\kappa;C(\kappa-4),\frac 12(\kappa-4))$ process started
from $(0;0^+,0^-)$. Let $K(\infty)=\cup_{t<\infty}K(t)$. Then $\pa_\HH^+ K(\infty)$ has the same distribution
as the image of a chordal SLE$(\kappa';(1-C)(\kappa'-4))$ trace started
from $(0;0^+)$, where $\kappa'=16/\kappa$.
\label{duality2*}
\end{Theorem}

\end{document}